\documentstyle[amssymb,amstex,a4wide,12pt]{article}
\newtheorem{theorem}{Theorem}
\newtheorem{lemma}{Lemma}
\newtheorem{corollary}{Corollary}

\newtheorem{remark}{Remark}

\newcommand{\Edges}{\operatorname{Edges}}

\newcommand{\card}{\operatorname{card}}
\newcommand{\conv}{\operatorname{conv}}
\newcommand{\lnaw}{\langle\langle}
\newcommand{\pnaw}{\rangle\rangle}

\newcommand{\lodc}{\overline{[[}}
\newcommand{\podc}{\overline{]]}}
\newcommand{\e}{\operatorname{e}}

\newcommand{\dist}{\operatorname{dist}}
\newcommand{\crop}{\operatorname{crop}}
\newcommand{\cropGraph}{\operatorname{CropGraph}}
\newcommand{\branch}{\operatorname{branch}}

\begin{document}

\setlength{\baselineskip}{1.3\baselineskip}

\author{Tomasz Schreiber\footnote{Mailing address: Tomasz Schreiber, Faculty of
 Mathematics \& Computer Science, Nicolaus Copernicus University,
 ul. Chopina 12 $\slash$ 18, 87-100 Toru\'n, Poland; tel.: (++48) (+56)
 6112951, fax: (++48) (+56) 6228979; e-mail: {\tt tomeks at mat.uni.torun.pl}},
 \footnote{
        Research supported by the Polish Minister of Science and Higher Education
        grant N N201 385234 (2008-2010)}\\
        Faculty of Mathematics \& Computer Science,\\
        Nicolaus Copernicus University, Toru\'n, Poland.}
        
\title{Polygonal web representation for higher order correlation functions 
       of consistent polygonal Markov fields in the plane}
\date{}
\maketitle

\paragraph{Abstract} 
 {\it We consider polygonal Markov fields originally introduced by Arak and Surgailis \cite{A1,AS1}.
      Our attention is focused on fields with nodes of order two, which can be regarded
      as continuum ensembles of non-intersecting contours in the plane, sharing a number
      of salient features with the two-dimensional Ising model. The purpose of this
      paper is to establish an explicit stochastic representation for the higher-order
      correlation functions of polygonal Markov fields in their consistency regime.
      The representation is given in terms of the so-called crop functionals (defined
      by a M\"obius-type formula) of polygonal webs which arise in a graphical
      construction dual to that giving rise to polygonal fields. The proof of
      our representation formula goes by constructing a martingale interpolation
      between the correlation functions of polygonal fields and crop functionals
      of polygonal webs.}


\section{Introduction}
   Polygonal Markov fields, originally introduced by Arak \& Surgailis \cite{A1,AS1}
   and further studied by Arak, Clifford \& Surgailis \cite{AS2,ACS} are random ensembles
   of non-intersecting polygonal contours in the plane interacting by hard-core exclusions
   and exhibiting two-dimensional germ-Markov property \cite{AS1}, with a variety
   of additional possible terms entering the Hamiltonian, including length
   and area elements (ibidem). The polygonal fields with V-shaped nodes (no nodes of order higher
   than two) as considered in this paper, share a number of essential features with
   the two-dimensional Ising model, prominent examples including the presence of 
   an Ising-like phase transition \cite{N1,SC05} as well as low temperature phase
   separation and Wulff droplet creation \cite{SC06}. For these reasons, the polygonal
   Markov fields are sometimes regarded as continuum counterparts of the Ising 
   model in the plane (as well as of the Potts model if higher order nodes are admitted).
   Remarkably, in many aspects the polygonal fields are {\it exactly tractable}, especially
   in the so-called {\it consistent regime} falling into the supercritical temperature region.
   In particular, at the consistency point one knows the exact value of the partition
   function, first and second order characteristics of the field \cite{AS1,AS2,ACS}.
   Further, a lot is known about the higher order correlations as well, including certain
   exact formulae \cite{SC08} and necessary and sufficient conditions for factorisation
   of edge correlations (ibidem) as well as an exponential mixing statement (asymptotic
   factorisation) for edge correlations (ibidem). A striking feature of polygonal Markov
   fields is that they admit a number of particularly convenient algorithmic constructions
   \--- graphical representations \cite{AS1,AS2,ACS,LS07,SC05,SC06,SC08}
   which are in fact the main tool for establishing of the afore-mentioned results.
   The geometric ingredient in these considerations is so predominant that often
   no supplementary calculations are needed, which stands in a strong contrast
   to the classical Ising model.
  
   The rich class of graphical constructions developed for polygonal fields have also
   found their applications in Bayesian image processing where they are used to generate
   image segmentations, see  \cite{CN94,KLS05,KLS07,LS07,SL08}. Experimenting with
   various black-white and grayscale images we already obtained promising results,
   further algorithms are a subject of our ongoing research in progress.

   The purpose of the present paper is to complement the existing exact results for consistent
   polygonal Markov fields by establishing in Theorem \ref{PWREPR} an explicit stochastic
   representation  for their higher order edge correlations in terms of expectations of
   the so-called crop functionals of {\it polygonal webs}. The polygonal web arises as the
   union of a collection of trajectories of continuous time critical branching polygonal
   random walks in the plane, interacting by a {\it bridge-creating} mechanism attempting
   to {\it clasp} random trees generated by branching walks into a {\it web} by establishing
   linear {\it bridges} between trees. As mentioned above, the polygonal fields admit
   a graphical construction whose full details in Section \ref{GDR} below.
   The dynamics of this construction can be interpreted as a (quite untypical)
   interacting particle system evolving in time and, speaking in this vein, the polygonal
   web arises in a graphical construction, given in Section \ref{PWEB} below, which can
   be regarded as {\it dual} to that for polygonal fields. The nature of this duality consists,
   roughly speaking, in the fact that the dynamic representation of polygonal webs exhibits
   features strongly reminiscent of those of the polygonal field construction under
   inverted time flow direction. To some extent this can be perceived as an analogy
   to the classical duality for interacting particle systems, see e.g. Section III.4 in
   \cite{LG}, although this is not a close analogy, expecially that here we deal with
   entire histories of the considered interacting particle systems (whose trajectories
   trace the polygonal fields and webs) rather than just with their instantaneous
   configurations. The proof of our duality representation for edge correlations
   goes by constructing a martingale interpolating between edge correlation
   functions and crop functionals of the corresponding polygonal webs.   
      
   The rest of this paper is organised as follows. In Section 
   \ref{NHPMF} below we present the concept of a non-homogeneous
   polygonal field as introduced in \cite{SC08}. Next, in Section
   \ref{GDR} we discuss the generalised graphical construction
   of such fields as developed ibidem. In the further Section
   \ref{PWEB} we develop a dual graphical construction and define
   the polygonal web. In the next Section \ref{MRES} we define
   edge correlation functions of polygonal fields and state
   our main representation Theorem \ref{PWREPR}
   as briefly discussed above. The proof of this theorem is
   given in Section \ref{PMT}, where a number of auxiliary
   constructions are also developed and many auxiliary results
   established. Finally, to the last Section \ref{PNCOR} we delagate
   the proof of a technical existence result for edge correlations.
      
\section{Non-homogeneous polygonal Markov fields in the plane and their consistent regime}\label{NHPMF}
 For an open bounded convex set $D$ define the family $\Gamma_D$ of admissible polygonal
 configurations on $D$ by taking all the finite planar graphs $\gamma$ in $D \cup \partial D,$
 with straight-line segments as edges, such that
 \begin{itemize}
  \item the edges of $\gamma$ do not intersect,
  \item all the interior vertices of $\gamma$ (lying in $D$) are of degree $2,$
  \item all the boundary vertices of $\gamma$ (lying in $\partial D$) are of degree $1,$
  \item no two edges of $\gamma$ are colinear.
 \end{itemize}
 In other words, $\gamma$ consists of a finite number of disjoint polygons, possibly nested
 and chopped off by the boundary. Further, for a finite collection $(l) = (l_i)_{i=1}^n$ of
 straight lines intersecting $D$ we write $\Gamma_D(l)$ to denote the family of admissible
 configurations $\gamma$ with the additional properties that
 $\gamma \subseteq \bigcup_{i=1}^n l_i$ and $\gamma \cap l_i$ is a single interval
 of a strictly positive length for each $l_i, i=1,...,n,$ possibly with some
 isolated points added.

 For a Borel subset of $A \subseteq {\Bbb R}^2$ by $[[A]]$ we shall denote the
 family of all straight lines hitting $A$ so that in particular $[[{\Bbb R}^2]]$
 stands for the collection of all straight lines in ${\Bbb R}^2.$ We shall also
 write $\lodc A \podc$ for the family of all linear segments in
 ${\Bbb R}^2$ hitting $A.$ Further, let $\mu$ be the standard isometry-invariant
 Haar-Lebesgue measure on the space $[[{\Bbb R}^2]]$ of straight lines in ${\Bbb R}^2.$
 Recall that one possible construction of $\mu$ goes by identifying a straight
 line $l$ with the pair $(\phi,\rho) \in [0,\pi) \times {\Bbb R},$
 where $(\rho \sin(\phi), \rho \cos(\phi))$ is the vector orthogonal to $l,$
 and joining it to the origin, and then by endowing the parameter space
 $[0,\pi) \times {\Bbb R}$ with the usual Lebesgue measure. Note that
 the above parametrisation of $[[{\Bbb R}^2]]$ with $[0,\pi) \times {\Bbb R}$
 endows $[[{\Bbb R}^2]]$ with a natural metric, topology and Borel $\sigma$-field
 which will be  used in this paper. Next, consider a non-negative Borel measure
 ${\cal M}$ on $[[{\Bbb R}^2]]$ admitting a locally bounded density $m(\cdot)$ with
 respect to $\mu.$ Below, the measure ${\cal M}$ will be interpreted as the activity
 measure on $[[{\Bbb R}^2]].$ Let $\Lambda^{\cal M}$ be the Poisson line process on
 $[[{\Bbb R}^2]]$ with intensity measure ${\cal M}$ and write $\Lambda^{\cal M}_D$
 for its restriction to the domain $D.$ Further, define the Hamiltonian 
 $L^{\cal M} : \Gamma_D \to {\Bbb R}_+$ given by 
 \begin{equation}\label{HAMILT}
  L^{\cal M}(\gamma) := \sum_{e \in \Edges(\gamma)} {\cal M}([[e]]),\;\; \gamma \in \Gamma_D.
 \end{equation}
 We note that the energy function $L^{\cal M}$ should be regarded as
 an anisotropic  environment-specific
 version of the length functional. Indeed, interpreting the activity ${\cal M}(dl)$ of
 a line $l$ hitting an edge $e \in \Edges(\gamma)$ as the likelihood of a new edge being
 created along $l$ intersecting and hence fracturing the edge $e$ in $\gamma,$
 we observe that, roughly speaking, the value of ${\cal M}([[e]])$ determines how likely
 the edge $e$ is to be fractured by another edge present in the environment.
 In other words, $L^{\cal M}(\gamma)$ determines {\it how difficult it is to
 maintain} the whole graph $\gamma \in \Gamma_D$ without fractures in the environment
 ${\cal M}$ \--- note that due to the anisotropy of the environment there may be graphs
 of a higher (lower) total edge length than $\gamma$ and yet of lower (higher) energy
 and thus easier (more difficult) to maintain and to keep unfractured due to the lack
 (presence) of high activity lines likely to fracture their edges.
   
 Following \cite{SC08}, with the above notation, for $\beta \in {\Bbb R}$ further
 referred to as the inverse temperature (from mathematical viewpoint also the unphysical
 negative values of inverse temperatures are admissible),
 we define the polygonal field ${\cal A}^{{\cal M};\beta}_D$ in $D$ with activity
 measure ${\cal M}$ to be the Gibbsian modification of the process induced on
 $\Gamma_D$ by $\Lambda^{\cal M}_D,$ with the Hamiltonian $L^{\cal M}$ at 
 inverse temperature $\beta,$ that is to say
 \begin{equation}\label{GIBBS1}
  {\Bbb P}\left( {\cal A}^{{\cal M};\beta}_D \in G \right)
  := \frac{{\Bbb E} \sum_{\gamma \in \Gamma_D(\Lambda^{\cal M}_D) \cap G}
     \exp\left(- \beta L^{\cal M}(\gamma) \right)}
          {{\Bbb E} \sum_{\gamma \in \Gamma_D(\Lambda^{\cal M}_D)}
     \exp\left(- \beta L^{\cal M}(\gamma) \right)}
 \end{equation}
 for all sets $G \subseteq \Gamma_D$ Borel measurable with respect to, say, the usual
 Hausdorff distance topology. Note that this definition can be rewritten as
 \begin{equation}\label{GIBBS2}
    {\Bbb P}({\cal A}^{{\cal M};\beta}_D \in d\gamma) \propto 
    \exp(- \beta L^{\cal M}(\gamma)) \prod_{e \in \Edges(\gamma)}
    {\cal M}(dl[e]),\; \gamma \in \Gamma_D,
 \end{equation}
 where $l[e]$ is the straight line extending $e.$ In other words, the probability
 of having ${\cal A}^{{\cal M};\beta}_D \in d\gamma$ is proportional to 
 the Boltzmann factor $\exp(-\beta L^{\cal M}(\gamma))$
 times the product of edge activities ${\cal M}(dl[e]),\; e \in \Edges(\gamma).$
 Observe also that this construction should be regarded as a specific version
 of the general polygonal model given in (2.11) of \cite{AS1}.
 The finiteness of the partition function
 \begin{equation}\label{PAFA}
   {\cal Z}^{{\cal M};\beta}_D := {\Bbb E} \sum_{\gamma \in \Gamma_D(\Lambda^{\cal M}_D)}
      \exp\left(- \beta L^{\cal M}(\gamma) \right) < \infty
 \end{equation}
 for all $\beta \in {\Bbb R}$ is not difficult to verify and has been established
 in \cite{SC08}, see (32) there.

 In this paper we shall focus on polygonal fields in their {\it consistent regime}
 corresponding to $\beta = 1.$ As shown in Section 3 in \cite{SC08}, this
 particular choice of temperature parameter places us in the context of a
 non-homogeneous version of Arak-Surgailis \cite{AS1} construction for the
 {\it consistent} polygonal fields, see Section 4 there. This ensures striking
 properties of the field. First of these, {\it the consistency}, states that
 for each open bounded and convex $D \subseteq D' \subset {\Bbb R}^2$ the field
 ${\cal A}^{{\cal M}}_D  := {\cal A}^{{\cal M};1}_D$ coincides in law with
 ${\cal A}^{{\cal M}}_{D'} := {\cal A}^{{\cal M};1}_{D'} \cap D,$
 thus allowing for a direct construction of the infinite volume process
 (thermodynamic limit) ${\cal A}^{\cal M} := {\cal A}^{{\cal M};1}$
 on the whole ${\Bbb R}^2$ such that ${\cal A}^{\cal M}_D = {\cal A}^{\cal M}
 \cap D.$ 
 %Note that this statement holds {\it only} for convex $D$ and $D'$
 %and, in general, it is false for non-convex domains!
 The infinite volume
 process ${\cal A}^{\cal M}$ takes its values in the space
 $\Gamma := \Gamma_{{\Bbb R}^2}$ of whole-plane admissible configurations,
 with obvious meaning of this notation. Further, the explicit
 formula for the partition function ${\cal Z}^{\cal M}_D$ is known for convex
 $D,$ see Theorem 4.1 in \cite{AS1} for the homogeneous case and Theorem 1 in \cite{SC08}
 for the general non-homogeneous set-up. We state this formula in (\ref{Zexp}).
 Moreover, one-dimensional linear sections of the field are fully
 characterised in distribution, see ibidem and Theorem 2 in \cite{SC08}.
 Finally, the polygonal fields ${\cal A}^{\cal M}_D$ enjoy the two-dimensional
 Markov property stating that  the conditional behaviour of the field inside a
 smooth closed curve depends on the outside configuration only through arbitrarily
 small neighbourhoods of the curve or, equivalently, through the trace of the
 external configuration on the curve, consisting of intersection points and
 directions. This property is a direct consequence of the Gibbsian definition
 (\ref{GIBBS1},\ref{GIBBS2}) of the field and, unlike the previous properties,
 it holds for all inverse temperatures $\beta \in {\Bbb R}$ rather than just
 for $\beta = 1.$ We will not discuss this property in the present paper
 and we refer the reader to the original work of Arak and Surgailis \cite{AS1}
 for further details.  

 \section{Generalised dynamic representation for consistent fields}\label{GDR}
 The present section is meant to recall the {\it generalised dynamic representation}
 for consistent polygonal fields as developed in Section 4 of \cite{SC08}, which will serve
 as a crucial tool in our further considerations. The name {\it generalised representation}
 comes from the fact that it generalises the original construction of such fields 
 introduced by Arak and Surgailis in \cite{AS1}. In the sequel we will often omit the
 qualifier {\it generalised} for the sake of terminological brevity. To describe
 the representation,
 fix the convex field domain $D$ and let $(D_t)_{t \in [0,1]}$ be a time-indexed
 increasing family of compact convex subsets of $\bar D,$ eventually covering the
 entire $\bar D$ and interpreted as a {\it growing window} gradually revealing increasing
 portions of the polygonal field under construction in the course of the time flow.
 In other words, under this interpretation, the portion of a polygonal field in a
 bounded open convex domain $D$ {\it uncovered} by time $t$ is precisely its
 intersection with $D_t.$ To put it in formal terms, consider $(D_t)_{t \in [0,1]}$
 satisfying
 \begin{description}
  \item{\bf (D1)} $(D_t)_{t \in [0,1]}$ is a strictly increasing family of compact
        convex subsets of $\bar D = D \cup \partial D.$
  \item{\bf (D2)} $D_0$ is a single point $x$ in $\bar{D} = D \cup \partial D.$
  \item{\bf (D3)} $D_1$ coincides with $\bar D.$
  \item{\bf (D4)} $D_t$ is continuous in the usual Hausdorff metric on compacts.  
 \end{description}
 Note that the extra fifth condition imposed on $D_t$ in Section 4 of \cite{SC08}
 is automatically satisfied here due to the absolute continuity ${\cal M} \ll \mu$ and
 thus is not mentioned here. Clearly, under these conditions, for ${\cal M}$-almost
 each $l \in [[D]]$ the intersection $l \cap D_{\tau_l}$ consists of precisely one
 point ${\Bbb A}(l),$ where $\tau_l = \inf \{ t \in [0,1],\; D_t \cap l \neq \emptyset \}.$ 
 The point ${\Bbb A}(l)$ will be referred to as the {\it anchor point} for $l,$
 this induces the {\it anchor} mapping ${\Bbb A}: [[D]] \to D$ defined ${\cal M}$-almost
 everywhere.
 Consider now the following dynamics in time $t \in [0,1],$ with all updates
 given by the rules below performed independently of each other.
 \begin{description}
  \item{\bf (GE:Initialise)} Begin with empty field at time $0,$
  \item{\bf (GE:Unfold)} 
                   Between critical moments listed below, during the time interval $[t,t+dt]$ the
                   unfolding field edges in $D_t$ reaching $\partial D_t$ extend straight to
                   $D_{t+dt} \setminus D_t,$
  \item{\bf (GE:BoundaryHit)}
                   When a field edge hits the boundary $\partial D,$ it stops growing in this direction
                   (note that ${\cal M}$-almost everywhere the intersection of a line with 
                    $\partial D$ consists of at most two points),
  \item{\bf (GE:Collision)}
                   When two unfolding field edges intersect in $D_{t+dt} \setminus D_t,$ they are not 
                   extended any further beyond the intersection point (stop growing in the 
                   direction marked by the interesction point),
  \item{\bf (GE:DirectionalUpdate)}
                   A field edge extending along $l \in [[D_t]]$ updates its direction during
                   $[t,t+dt]$ and starts unfolding along $l' \in [[l^{[t,t+dt]}]],$
                   extending away from the anchor point ${\Bbb A}(l'),$ with probability
                   ${\cal M}(dl'),$ where $l^{[t,t+dt]} := l \cap (D_{t+dt} \setminus D_t).$
                   Directional updates of this type are all performed independently,
  \item{\bf (GE:LineBirth)} Whenever the anchor point ${\Bbb A}(l)$ of a line $l$ falls
                  into $D_{t+dt} \setminus D_t,$ the line $l$ is born at the time $t$ at
                  its anchor point with probability ${\cal M}(dl),$ whereupon it begins
                  extending in both directions with the growth of $D_t$ (recall that
                  $l$ is ${\cal M}$-almost always tangential to $\partial D_t$ here),
  \item{\bf (GE:VertexBirth)} For each intersection point of lines $l_1$ and $l_2$ falling
                  into $D_{t+dt} \setminus D_t,$ the pair of field lines $l_1$ and $l_2$ is
                  born at $l_1 \cap l_2$ with probability ${\cal M}(dl_1) {\cal M}(dl_2),$
                  whereupon both lines begin unfolding in the directions away from their
                  respective anchor points ${\Bbb A}(l_1)$ and ${\Bbb A}(l_2).$ 
 \end{description}
 Observe that the evolution rule {\bf (GE:VertexBirth)} means that pairs of
 lines are born at birth sites distributed according to a Poisson point process in
 $D$ with intensity measure given by the {\it intersection measure} $\lnaw {\cal M} \pnaw$
 of ${\cal M}$ defined as follows 
 \begin{equation}\label{MM}
  \lnaw {\cal M} \pnaw(A) := 
  \frac{1}{2} {\cal M} \times {\cal M}( \{ (l_1,l_2),\; l_1 \cap l_2 \subset A \}). 
 \end{equation}
 The importance of the intersection measure lies in the fact that
 \begin{equation}\label{Zexp}
  {\cal Z}^{{\cal M};1}_D = \exp(\lnaw {\cal M} \pnaw(D))
 \end{equation} 
 as shown in Theorem 1 in \cite{SC08}.
 The following result stating that the polygonal field resulting from the above
 construction actually coincides with ${\cal A}^{\cal M}_D$ has been established
 in \cite{SC08}, see Theorem 3 there. 
 \begin{theorem}\label{GEthm}
  The random contour ensemble resulting from the above construction {\bf (GE)}
  coincides in law with ${\cal A}^{\cal M}_D.$
 \end{theorem}

\section{Polygonal web}\label{PWEB}
 Having defined the non-homogeneous polygonal fields and presented their graphical
 construction, we pass now to another object central to this paper, which we
 name the {\it polygonal web}. Whereas the details of the connection between the critical 
 polygonal web and the corresponding consistent polygonal field are to be established
 in the subsequent Section \ref{MRES}, here we emphasise that, in a sense, the polygonal web
 constitutes the {\it dual object} to the polygonal field sharing the same activity
 measure, and this duality is going to be reflected in the construction of the
 polygonal web. Roughly speaking, the polygonal web arises as the union of 
 critical branching polygonal random walks, interacting by an additional
 {\it bridge-creating} mechanism, clasping the branched polygonal trees into a web.

 \subsection{Constructing the polygonal web}
  Consider an open bounded and convex domain $D$ and let $(D_t)_{t \in [0,1]}$ be
  a growing family of compact subsets of $\bar D$ satisfying the usual conditions {\bf (D1-4)}
  as in Section \ref{GDR}. Further, assume a collection $(\bar l,\bar x) = (l_i,x_i)_{i=1}^k$
  is given where $x_i$ are points in $\bar D$ whereas $l_i$ are straight lines with
  $x_i \in l_i,\; i=1,\ldots,k.$ Moreover, to avoid uninteresting pathologies we require
  that no three different lines $l_i,\;i=1,\ldots,k,$ intersect at one point.
  Each such pair $(l_i,x_i)$ will be called 
  an {\it edge marker} because $l_i$ can be interpreted as a directional indicator
  for a linear segment/edge passing through $x_i,$ see Section \ref{MRES} below
  where this interpretation is further developed and exploited. The entire collection $(\bar l,\bar x)$
  will be referred to as the {\it edge marker configuration}. The polygonal web
  ${\cal W}[(\bar l,\bar x)] := {\cal W}^{\cal M}_D[(\bar l,\bar x)]$
  generated by $(\bar l,\bar x)$ in
  $\bar D,$ with activity measure ${\cal M},$ is the union of polygonal trees
  in $\bar D$ arising as the final state $w_1$ of the following graphical
  construction process $w_s$ evolving for $s \in [0,1],$ where all random 
  updates listed are performed independently.
  \begin{description}
   \item{\bf (W:Start)}
        At the time $s = 0$ we let $w_0$ consist of zero-length edges (edge germs)
        at $x_i,\; i=1,\ldots,k,$ directed along the respective $l_i$'s.
   \item{\bf (W:GrowInwards)}
        Between the critical moments listed below, during the time interval $[s,s+ds]$
        all edges of $w_s$ reaching the boundary $\partial D_{1-s}$ extend straight
        to $D_{1-s} \setminus D_{1-s-ds}$ along their respective directional lines.
        The edges (edge germs) not yet touched by the
        boundary of the shrinking domain $D_{1-s}$ (and hence contained in the
        interior of $D_{1-s}$) remain intact and do not evolve
        until eventually hit by the boundary at some later time (unless they get
        {\it frozen} prior to that, see below, in which case they never start evolving).
        Below, we call
        edges reaching the current boundary $\partial D_{1-s}$ {\it active} and
        we say that edges (edge germs) not yet hit by the boundary are {\it inactive}.
        Note that the intersection point of a web edge with the current boundary
        $\partial D_{1-s}$ can be interpreted as its instantaneous {\it growth-tip}
        and so will it be called in the sequel. Observe that inactive edges do not
        have growth-tips.   
  \item{\bf (W:BranchAndTurn)}
        During the time interval $[s,s+ds]$ an active web edge reaching the boundary
        $\partial D_{1-s}$ and extending along $l \in [[D_{1-s}]]$ yields a new offspring
        edge starting at $l^{[s,s+ds]} := l \cap D_{1-s} \setminus D_{1-s-ds}$
        %$l \cap \partial D_{1-s}$ 
        and directed along
        $l' \in [[l^{[s,s+ds]}]]$ %,\; l^{[s,s+ds]} := l \cap D_{1-s} \setminus D_{1-s-ds},$
        with probability ${\cal M}(dl').$ Both the original and offspring edges go on
        evolving according to the usual rules. The branching updates are performed
        independently for all active web edges.
  \item{\bf (W:ForcedBranchAndTurn)}
        If during the time interval $[s,s+ds]$ an active web edge extending along
        some $l \in [[D_{1-s}]]$ and reaching the  boundary $\partial D_{1-s}$
        intersects the directional line $l'$ of some other non-frozen web edge in
        $w_s \cap D_{1-s},$ be it active or inactive, and there is no web edge
        along $l'$ reaching $l \cap \partial D_{1-s}$ and created in a prior
        stage of the evolution, then a new offspring edge
        is created at $l \cap \partial D_{1-s}$ directed along $l'$ and both the
        original  and offspring edges go on evolving according to the usual rules.
        Note that we only consider directional lines $l'$ of edges hitting the
        domain $D_{1-s}$ \--- the edges present in $\bar D \setminus D_{1-s}$ but
        terminated before the time $s$ (see below for edge termination events) are 
        not taken into account. Thus, the forced branching occurs if the growth-tip
        of an edge hits the current directional line of another (non-frozen) edge
        currently present in
        the system. Whenever a forced branching occurs, the resulting
        offspring edge is called a {\it forced edge} whereas the edge whose directional
        line $l'$ gave rise to the branching is referred to as the corresponding {\it forcing
        edge}. Observe that a pair of forced and forcing edges will meet and
        coalesce into a single edge at further stages of the construction unless
        one of the edges terminates prior to that.

        Note that if the forcing edge is active, it resides at the boundary of
        the domain $D_s$ and yields a single forced edge on the opposite side of
        the domain. On the other hand, an inactive forcing edge (edge germ)
        located in the interior of $D_s$ may give rise to (at most) two forced
        edges, one on each side of the domain.
  \item{\bf (W:Terminate)}
        During  the time interval $[s,s+ds]$ an active web edge reaching the boundary
        $\partial D_{1-s}$ and extending along $l \in [[D_{1-s}]]$ terminates (stops
        evolving) with probability ${\cal M}([[l^{[s,s+ds]}]])$ where 
        $l^{[s,s+ds]} := l \cap D_{1-s} \setminus D_{1-s-ds}.$
  \item{\bf (W:StopIfSeparated)}
        Whenever at some time moment $s \in [0,1]$ a web edge (or edge germ) $e$
        in $w_s \cap D_{1-s},$ be it active or inactive, has the property that
        $l[e] \cap \conv([w_s \setminus e] \cap \partial D_{1-s}) = \emptyset$
        (that is to say the directional line $l[e]$ of the web edge $e$ does not
         hit the convex hull generated by the current growth-tips of the remaining
         non-terminated 
         active edges and non-frozen germination points of the remaining inactive edges, 
         in which case we say that $e$ {\it separates} from $w_s$ at the time $s$),
         then $e$ terminates and stops evolving at this
         point. Note that in case of an inactive edge germ $e$ {\it to terminate} means
         to remain {\it frozen} in inactive state and never to activate even when hit
         by the boundary at the later stages of the evolution.
 \end{description}
  A careful reader might ask at this moment why in the above construction we do not 
  consider the case when at some time $s$ a web edge becomes tangential to the
  domain $D_{1-s}.$ The answer is that, with probability one, such cases do
  not occur in the course of the construction because an edge to become
  tangential to the boundary of the domain separates from the web prior to
  that and thus gets terminated by an application of {\bf (W:StopIfSeparated)}
  rule.   
   
  The construction of the polygonal web as presented above admits a natural 
  intuitive description \--- the edge germs initiating the process emit 
  critical branching polygonal random walks directed by the activity measure
  ${\cal M}$ and unfolding {\it inward} the domain $D$ which can be regarded
  as dual to the dynamic representation in Section \ref{GDR} where the growth
  was directed {\it outwards}. The branching is critical because the binary
  branching and termination intensities coincide in {\bf (W:BranchAndTurn)} and
  {\bf (W:Terminate)}. The role of the additional {\bf (W:ForcedBranchAndTurn)}
  rule is to ensure the possibility of {\it bridging} the gaps between two
  separated co-linear
  parts of the same segment present on two opposite sides of the window $D_{1-s},$
  thus clasping the polygonal trees into a web. Finally, as may 
  be seen in the sequel, the {\bf (W:StopIfSeparated)} rule reflects the
  structural knowledge on independence of edge covering events in
  polygonal fields as established in \cite{SC08}, and as such it
  is not indispensible in its full form for the theory developed
  below to be valid and may be replaced by various weaker variants,
  see Remark \ref{Sepa}.
  
  It is useful to note that, regarded as a polygonal graph, the polygonal
  web contains T-shaped nodes (branching points), I-shaped nodes (edge terminal points)
  and X-shaped nodes (edge intersection points) but no V-shaped nodes.

 \subsection{Crop functional}
  To establish a direct link between the polygonal web and polygonal fields
  we define now the {\it crop functional} of the polygonal web, further denoted
  as $\crop({\cal W}[(\bar l,\bar x)]).$ To this end, we identify the polygonal
  web ${\cal W}[(\bar l,\bar x)] = {\cal W}^{\cal M}_D[(\bar l,\bar x)]$
  with the collection of {\it branches} connecting
  the initial edge germ locations $x_i,\; i=1,\ldots,k$ to the edge termination
  points $y_1,\ldots,y_m$ resulting from {\bf (W:Terminate,StopIfSeparated)} rules or
  arising as meeting points where forcer-forced pairs of co-linear web edges
  unfolding in opposite directions merge into linear segments. To avoid nuisance
  technicalities below we formally interpret each meeting point of a forcer-forced
  pair as {\it two distinct points}, one terminating the forcing branch, the second
  terminating the forced branch. Clearly, $m \geq k$ since each initial
  edge germ $x_i$ emits at least one branch and each branch eventually
  terminates. Keeping $(\bar l,\bar x)$ fixed, we shall index the branches
  constituting ${\cal W}[(\bar l,\bar x)]$ by their terminal points,
  writing $\branch[y_j]$ for the branch terminating at $y_j$ \--- the
  inambiguity of this indexation is ensured by the above convention
  on meeting points of forcer-forced pairs. Along
  each branch we have a natural {\it chronological} ordering from the {\it root}
  $x_i,\; i \in \{1,\ldots,k\}$ to the {\it endpoint} $y_j,\; j \in \{1,\ldots,m\}.$
  For a collection of branches $\branch[y],\; y \in {\cal Y},\; {\cal Y} \subseteq
  \{ y_1,\ldots,y_m \},$ we consider the induced  polygonal {\it crop graph} 
  $\cropGraph[{\cal Y}]$ obtained as follows.
  \begin{description}
   \item{\bf (Crop:Grow)} 
         Follow the growth of all branches in $\{ \branch[y],\; y \in {\cal Y}\}$
         starting from their roots and unfolding towards their respective endpoints
         during the time interval $[0,1]$ as in the course of the polygonal
         web dynamics {\bf (W)}.
   %\item{\bf (Crop:StopIfSeparated)}
   %      Whenever in the course of the dynamics a branch (including an inactive edge germ)
   %      separates from the remaining non-terminated branches, it stops growing at this point.
   \item{\bf (Crop:StopOnCollision)}
         Whenever in the course of their growth two branches meet, they both stop
         growing at this point.
  \end{description}
  Thus, the crop graph is a subgraph of ${\cal W}[(\bar l,\bar x)]$ containing
  T-shaped, I-shaped and V-shaped nodes but not X-shaped nodes. This is because
  the crop graph arises by (recursively) cutting off the parts of branches past
  their intersections with other branches present in the inducing collection.
  There are two ways in which two branches can meet in the  above construction
  \--- they can either intersect coming from two non-colinear directions and
  giving rise to a V-shaped node in $\cropGraph[{\cal Y}],$ or meet coming
  from opposite co-linear directions yielding a linear segment rather than
  a graph node. Note that two distinct branches sharing a common sub-branch
  and thus coinciding during initial growth phase are not considered to meet
  or intersect! On the other hand, if several distinct branches coinciding 
  during an initial growth phase intersect another branch(es) during this
  phase, the growth-interrupting {\bf (Crop:StopOnCollision)} rule applies to
  {\it all} these branches {\it simultaneously}.% Finally, note also that
  %the {\it separation} of a branch from the remaining ones, as considered
  %in {\bf (Crop:StopIfSeparated)} here, is understood
  %in full analogy to {\bf (W:StopIfSeparated)} above
  %\--- at some time $s$ the current directional line of one branch fails to hit
  %the convex hull generated by the current growth-tips of the remaining 
  %active branches (not coinciding with the one considered) and non-frozen 
  %germination points of the remaining inactive edges.  

  %%meeting points of the remaining growing 
  %% branches with $\partial D_{1-s}$ or equivalently $D_{1-s}.$

  It is easily seen that two different collections of terminal points can yield
  identical crop graphs because some branch turning points can be cut off due to
  collisions in
  the course of the {\bf (Crop)} dynamics. It is clear though that for each
  instance of a crop graph arising in {\bf (Crop)} dynamics there exists a
  unique collection of terminal points with the property that no branch
  turning points occur past the cut-off points. This unique collection is
  called {\it minimal} for its crop graph, or just {\it minimal} for short if
  no ambiguity arises.  

  To proceed, we say that a collection ${\cal Y} \subseteq \{ y_1,\ldots,y_m \}$
  of branch-determining endpoints is {\it complete} iff each initial germ location
  is the root of some $\branch[y],\; y \in {\cal Y}.$ Clearly, the cardinality
  of such a complete collection cannot fall below $k.$ Further, we say that the
  crop graph $\cropGraph[{\cal Y}]$ of a complete endpoint collection ${\cal Y}$
  is {\it normal} iff it contains no forced edges which fail to eventually meet
  and merge with their co-linear forcing edges. 
  %
  %\begin{itemize}
  % \item The collection ${\cal Y}$ is minimal for $\cropGraph[{\cal Y}].$
  % \item In addition,  the graph contains no forced edges which fail to
  %       eventually meet and merge with their co-linear enforcing edges.
  %\end{itemize}
 Note that this condition can be violated by either having a forced edge without its
 forcer present in the graph or due to a death or directional update along either a
 forcing or a forced edge. Whereas the death of either edge in an forcer-forced pair
 does necessarily lead to the lack of normality, a directional update does so only if there
 is no other branch in the collection along which the considered forcer-forced pair could
 extend further past the turning point. A graph which is not normal is called {\it abnormal}.
 Thus, roughly speaking, an abnormal graph is a graph containing the forced part of
 an {\it incomplete bridge}. Write 
  \begin{equation}\label{IOTA}
    \iota({\cal Y}) := \left\{ \begin{array}{ll} 1, & \mbox{ if ${\cal Y} \subseteq \{y_1,\ldots,y_m\}$
    is complete, minimal and $\cropGraph[{\cal Y}]$ is normal, } \\
    0, & \mbox{ otherwise. } \end{array} \right.
  \end{equation} 
  With this notation, we define the crop functional
  \begin{equation}\label{CROPFN}
   \crop({\cal W}[(\bar l,\bar x)]) :=
   \sum_{{\cal Y} \subseteq \{ y_1,\ldots,y_m \}} 
   (-1)^{\card({\cal Y})-k} \iota({\cal Y}).
  \end{equation}
  The expression (\ref{CROPFN}) has the aesthetic advantage of
  defining the crop functional in a form reminiscent of the classical
  inverse M\"obius transform, with the summation performed over all subsets
  of $\{ y_1,\ldots,y_m \},$ see e.g. Section 2.6 in \cite{MM}. To exploit
  this feature define
  $$ \hat{\iota}({\cal Y})  
    := \left\{ \begin{array}{ll} 1, & \mbox{ if $\cropGraph[{\cal Y}]$ is normal, } \\
    0, & \mbox{ otherwise, } \end{array} \right. 
  $$
  that is to say $\hat{\iota}$ is the indicator of crop graph normality, without the
  extra completeness and minimality requirements. Given a crop graph $\varrho
  = \cropGraph[{\cal Y}]$ for some complete and minimal ${\cal Y}_{\varrho} \subseteq \{y_1,\ldots,y_m\},$ 
  we let $B[y],\; y \in {\cal Y}_{\varrho},$ be the set of all $y_j,\;j=1,\ldots,m,$ such that
  $\branch[y_j]$ contains the entire subbranch $\branch[y] \cap \varrho.$
  Then, for any ${\cal Y} \subseteq \{y_1,\ldots,y_m\}$
  such that $\cropGraph[{\cal Y}] = \cropGraph[{\cal Y}_{\varrho}]$
  we have $\hat\iota({\cal Y}) = \iota({\cal Y}_{\varrho})$ and, moreover, ${\cal Y}$ decomposes into
  the disjoint union of non-empty ${\cal Y}_y = {\cal Y} \cap B[y],\; y \in {\cal Y}_{\varrho}.$ 
  We also have $\card({\cal Y}) = \card({\cal Y}_{\varrho}) + \sum_{y \in {\cal Y}_{\varrho}}
  [\card({\cal Y}_y) - 1].$
  Consequently, using Newton's binomial formula,
  $$ \sum_{{\cal Y},\; \cropGraph[{\cal Y}] = \cropGraph[{\cal Y}_{\varrho}]}
     (-1)^{\card({\cal Y})-k} \hat\iota({\cal Y}) = $$ $$ (-1)^{\card({\cal Y}_{\varrho})-k} 
     \iota({\cal Y}_{\varrho})
     \prod_{y \in {\cal Y}_{\varrho}} \sum_{\emptyset \neq {\cal Y}_y \subseteq B[y]}
     (-1)^{\card({\cal Y}_y)-1} =
     (-1)^{\card({\cal Y}_{\varrho})-k} \iota({\cal Y}_{\varrho}) \prod_{y \in {\cal Y}_{\varrho}} 1 = $$
 $$  (-1)^{\card({\cal Y}_{\varrho})-k} \iota({\cal Y}_{\varrho}).
 $$
  Thus, (\ref{CROPFN}) can be alternatively rewritten as
  \begin{equation}\label{CROPFNPOM}
   \crop({\cal W}[(\bar l,\bar x)]) =
   \sum_{{\cal Y} \subseteq \{ y_1,\ldots,y_m \},\; {\cal Y} \; \mbox{ complete}} 
   (-1)^{\card({\cal Y})-k} \hat\iota({\cal Y}).
  \end{equation}     
  Further, we put
  $T[x_i] := \{ y_j \in \{ y_1,\ldots,y_m\},\; {\rm root}(\branch[y_j]) = x_i \},\; i=1,\ldots,k,$
  and write, applying Newton's binomial formula,
  $$ \sum_{{\cal Y} \subseteq \{ y_1,\ldots,y_m \},\; {\cal Y} \; \mbox{complete}}
     (-1)^{\card({\cal Y})-k} = \prod_{i=1}^k \sum_{\emptyset \neq {\cal Y}_i \subseteq T[x_i]}
     (-1)^{\card({\cal Y}_i)-1} = \prod_{i=1}^k 1 = 1. $$
  Combining this with (\ref{CROPFNPOM}) we conclude that
  \begin{equation}\label{CROPFN2}
   \crop({\cal W}[(\bar l,\bar x)]) = 1 - \sum_{{\cal Y} \subseteq \{ y_1,\ldots,y_m \},\;
   {\cal Y} \; \mbox{complete}}
   (-1)^{\card({\cal Y})-k} [\hat\iota({\cal Y})-1]
  \end{equation}
  and hence the crop equals {\it one} for polygonal webs whose all complete branch
  subcollections yield normal crop graphs, whereas the deviations of the crop
  functional from the {\it standard} value one are due to crop graph abnormalities,
  which will be further exploited in the sequel.

  A natural alternative way of defining the crop functional involves summation
  over crop subgraphs of ${\cal W}[(\bar l,\bar x)],$ that is to say over all
  possible {\it different} graphs arising in the {\bf (Crop)} dynamics above, in
  which case (\ref{CROPFN}) becomes
  \begin{equation}\label{CROPFN3}
   \crop({\cal W}[(\bar l,\bar x)]) = \sum_{\varrho \; \mbox{ is a normal crop graph
   in } {\cal W}[(\bar l,\bar x)]} (-1)^{\mbox{ number of branchings in } \varrho}
  \end{equation}
  since the number of branchings in $\varrho = \cropGraph[{\cal Y}],\; {\cal Y}$ minimal,
  is easily seen to coincide with $\card({\cal Y}) - k.$     

 \subsection{Interpretation of the crop functional}\label{ICF}
  A few words are due at this point to provide an intuitive interpretation of the
  crop functional as formally defined in (\ref{CROPFN}). To this end, we begin by
  mentioning a kid game quite popular in the happy time of the author's childhood:
  given a collection of arrows on a sheet of paper draw a family of closed curves passing
 through  these arrows, and in case where this can be done in more than one way resolve
  the ambiguity by trying to make the resulting picture resemble some real-life
  object. In fact, this and related problems are not just games and find serious
  interest in studies on human and computer vision, see \cite{AC07} and the
  references therein. In mathematical terms and specialising to our polygonal
  set-up, given an edge marker collection $(\bar l, \bar x) = (l_i,x_i)_{i=1}^k$
  we ask for  admissible polygonal configurations $\gamma \in \Gamma_{{\Bbb R}^2}$
  with the following properties
  \begin{itemize}
   \item Each edge $e$ of $\gamma$ contains some edge marker point $x_{i(e)}$
         such that $l[e] = l_{i(e)}.$
   \item Each edge marker point $x_i$ is contained in some $e(i) \in \Edges(\gamma)$
         such that $l[e(i)] = l_i.$
  \end{itemize}
  Following Section 5 of \cite{SC08} we denote the family of such configurations
  by $\Gamma(\bar l,\bar x)$ and write $N(\bar l,\bar x)$ for the cardinality
  of $\Gamma(\bar l,\bar x),$ that is to say the number of solutions to the
  discussed problem. To proceed, consider first a particularly simple
  deterministic instance ${\cal W}^0[(\bar l,\bar x)] := {\cal W}^0_D[(\bar l,\bar x)]$
  of polygonal web generated by $(\bar l,\bar x)$ not depending on the activity measure
  ${\cal M}$ \--- define ${\cal W}^0[(\bar l,\bar x)]$ to arise in the course of the
  above {\bf (W)} dynamics without non-forced turns/branchings and without termination
  events, which is a usual situation for example when $x_i$'s are very close to each
  other and are all contained in a  domain $D$ with very small ${\cal M}([[D]])$
  where applications of {\bf (W:TurnAndBranch,Terminate)} are very unlikely. Observe
  that the notation ${\cal W}^0_D[(\bar l,\bar x)]$ comes from the fact that the
  above construction of this polygonal web coincides with the {\bf (W)} dynamics
  under zero activity measure. Then, as will be shown in Lemma \ref{WN} below,
  $$ \crop({\cal W}^0_D[(\bar l,\bar x)]) = N(\bar l,\bar x) $$
  for each bounded convex domain $D$ containing all $x_i$'s.

  At the other extremity lies the typical behaviour of the crop
  functional of polygonal webs generated by marker configurations
  $(\bar l,\bar x)$ where the distances between all $x_i$'s are huge
  and where all $l_i$'s are pairwise different 
  \--- due to the criticality of our branching mechanism the usual
  situation then is that the web ${\cal W}[(\bar l,\bar x)]$
  splits into disjoint and distant sub-webs originating from respective
  $x_i$'s and the bridge creating attempts between these sub-webs fail
  with overwhelming probability in consequence of {\bf (W:Terminate)}. 
  In these circumstances all normal crop graphs in ${\cal W}[(\bar l,\bar x)]$ are readily seen
  to be unions of disjoint normal crop graphs in sub-webs stemming
  from individual $x_i$'s and thus, by (\ref{CROPFN3}), the value
  of $\crop({\cal W}[(\bar l,\bar x)])$ factorises into the
  product of crops of the sub-webs. Using (\ref{NISKIEK}) 
  below combined with our main Theorem \ref{PWREPR} and
  resorting to Remark \ref{Sepa} allows us to conclude easily
  but not immediately that the expectations of these individual
  crops are all $1$ and consequently
  $$ {\Bbb E}\crop({\cal W}[(\bar l,\bar x)]) \approx 1 $$
  as well.
  This fact, not studied in detail here, is intimately related to mixing
  properties of polygonal fields (asymptotic factorisation of edge correlations)
  but this link will not be followed any further in this paper because so far we
  are only able to establish a {\it slow polynomial} mixing using a rather technical 
  argument based on polygonal webs, whereas alternative methods
  developed in Section 7 of \cite{SC08} allowed us to establish
  {\it exponential} mixing there at least for rectangular fields. 
  
  In remaining situations the crop functional interpolates between the
  above extremities.  
      
\section{Polygonal web representation for edge correlations}\label{MRES}
 The purpose of this section is to formulate our main result stating that
 arbitrary order edge correlation functions of the polygonal field coincide
 with the expectations of the respective polygonal web crops. 

 \paragraph{Edge correlations}
 Having introduced the crucial concepts in preceding sections we are now in
 a position to define the principal object of our study in this paper,
 that is to say the edge-correlation functions for polygonal fields,
 and to formulate our main results. Due to the polygonal nature of the
 considered field the natural object to consider are the {\it edge correlations}
 \begin{equation}\label{CORR}
  \sigma^{\cal M}[dl_1,x_1;\ldots;dl_k,x_k] := {\Bbb P}\left( \forall_{i=1}^k 
  \exists_{e \in \Edges({\cal A}^{\cal M})} \; \pi_{l_i}(x_i) \in e,\; l[e] \in dl_i \right),
 \end{equation}
 where $l_1,\ldots,l_k$ are straight lines and $\pi_{l_i}$ is the orthogonal
 projection on $l_i.$ In all cases below we shall be interested in correlations
 with $x_i \in l_i,$ in which case $\sigma^{\cal M}[dl_1,x_1;\ldots;$ $dl_k,x_k]$ 
 can be interpreted as the probability element that the polygonal field ${\cal A}^{\cal M}$
 passes through points $x_i$ in the directions determined by the respective lines
 $l_i,\; i=1,\ldots,k.$ For general $x_i,$ not necessarily lying on $l_i,$ the $k$-fold
 correlation $\sigma^{\cal M}[dl_1,x_1;\ldots;dl_k,x_k]$ is the probability that
 the polygonal field passes through points $\pi_{l_i}(x_i)$ in the directions
 determined by the respective lines $l_i,\; i=1,\ldots,k.$ 

 Recall from Section \ref{PWEB} above that collections $(\bar l,\bar x) = (l_i,x_i)_{i=1}^k$
 of lines $l_i$ and points $x_i$ with $x_i \in l_i$ are referred to as
 {\it edge marker configurations (collections)}
 whereas each pair $(l_i,x_i)$ belonging to such a collection is called
 an {\it edge marker}. No edge marker can occur twice in an edge marker collection. 
 An edge along $l_i$ passing through $x_i$ is said to 
 {\it cover the marker} $(l_i,x_i).$
 We say that an edge marker collection $(\bar l,\bar x) =
 (x_i,l_i)_{i=1}^k$ is in general position if the lines $l_i$ are pairwise
 different and $x_j \not\in l_i$ for $j \neq i,$ otherwise if $l_i = l_j$
 for some $i \neq j$ then the collection is called degenerate and $x_i, x_j$
 are declared {\it coupled} by $l_i = l_j,$ finally
 if $x_i \in l_j$ for some $i \neq j$ with $l_i \neq l_j$ then the collection
 is said to be in singular position. Thus, an edge marker collection can be
 simultaneously degenerate and singular. As mentioned above, if two edge marker lines $l_i = l_j$
 in a degenerate configuration coincide, we say that the edge markers $(l_i,x_i)$
 and $(l_j,x_j)$ are {\it coupled}, sometimes for brevity we just say that
 $x_i$ and $x_j$ are coupled. While allowing both for singularity and degeneracy
 of marker collections, we strictly exclude the situations where three or more
 {\it different} marker lines intersect at one point, in order to avoid unnecessary
 technical pathologies. 

 Note that in the singular case
 where $x_i \in l_j$ for some $i \neq j$ but $l_i \neq l_j$ it makes sense to
 consider {\it one-sided} correlations  $\sigma^{\cal M}[dl^+_i,x_i;dl_j,x_j;\ldots]$
 and $\sigma^{\cal M}[dl^-_i,x_i;dl_j,x_j;\ldots]$ with $l_i^+$ and $l_i^-$
 standing for two half-lines into which $l_i$ is split by the intersection
 point $\{ x_i \} = l_i \cap l_j.$  The definition of such correlations
 is analogous to (\ref{CORR}) the difference being that the field edge
 containing $x_i$ is required to extend respectively at least in the
 direction of $l_i^+$ and $l_i^-,$ yet it is also allowed although not
 required to extend in the opposite direction as well. Since ${\cal M} \ll \mu,$ it
 is easily seen that for a configuration in general position we would have
 $\sigma^{\cal M}[dl^+_i,x_i;dl_j,x_j;\ldots]  = 
  \sigma^{\cal M}[dl^-_i,x_i;dl_j,x_j;\ldots] = \sigma^{\cal M}[dl_i,x_i;dl_j,x_j;\ldots],$
 but this is no more the case in the considered singular situation, where the respective 
 correlations are non-trivially affected by the event that an edge along $l_j$ may extend
 from $x_j$ down to $x_i$ where it may intersect with another edge along $l_i.$

 In the
 sequel we will use the notation $\Delta^{*}[k],\; * \in \{ g, s, d \},$ for the
 respective sets of all collections $(l_i,x_i)_{i=1}^k$ in general, singular
 and degenerate positions. Further, we
 put $\Delta^* := \bigcup_{k=1}^{\infty} \Delta^*[k],\;
 * \in \{ g,s,d \}$ and $\Delta[k] := \bigcup_{* \in \{ g,s,d \}} \Delta^*[k]$
 and $\Delta := \bigcup_{* \in \{ g,s,d \}} \Delta[k].$ 
 Clearly, $\Delta[k]$ can be endowed with the natural
 product topology from $([[{\Bbb R}^2]] \times {\Bbb R}^2)^{\times k}.$  

 \paragraph{Edge correlation functions}
 For a non-degenerate edge marker configuration 
 $(\bar l,\bar x) = (l_i,x_i)_{i=1}^k$ the (normalised) k-fold edge correlation
 function $\phi(l_1,x_1;\ldots;l_k,x_k)$ is defined by
 \begin{equation}\label{MDEFI}
  \phi(l_1,x_1;\ldots;l_k,x_k) = 
  \phi^{\cal M}(l_1,x_1;\ldots;l_k,x_k) := \frac{\sigma^{\cal M}[dl_1,x_1;\ldots;dl_k,x_k]}
  {{\cal M}(dl_1) \ldots {\cal M}(dl_k)}.
 \end{equation}  
 More generally, if $(\bar l,\bar x)$ is a degenerate edge marker configuration,
 its normalised correlation is given by
 \begin{equation}\label{MDEFI2}
  \phi(\bar l,\bar x) = \frac{\sigma^{\cal M}[dl_1,x_1;\ldots;dl_k,x_k]}{\prod^* {\cal M}(dl_i)}
 \end{equation}
 where $\prod^*$ stands for the product over lines $l_i$ in which each line is present exactly
 once, with repetitions discarded. We also adopt the convention that $\phi(\emptyset) := 1,$
 where $\emptyset$ stands for empty edge marker configuration.
 The existence of edge correlation functions is guaranteed
 by the following lemma.
 \begin{lemma}\label{NCOR}
  The correlation function $\phi(\bar l,\bar x)$ exists for 
  each edge marker configuration $(\bar x,\bar l)$ and is continuous
  on $\Delta^g[k].$ In particular, $\phi(\cdot)$ is locally bounded.
  % and satisfies
  % $\phi(\bar l,\bar x) \leq k!$ 
  % and is uniformly bounded
  % for all $(\bar l,\bar x) \in \Delta^g[k] \cup \Delta^s[k],\; k=1,2,\ldots$
 \end{lemma}
 Likewise, we can also consider one-sided versions of the correlation functions,
 for which the existence and local boundedness statements of Lemma \ref{NCOR}
 readily extend. Note that in general the correlation functions are not continuous
 at singular configurations in $\Delta^s[k]$ \--- indeed, if a sequence of configurations
 in general positions converges to some singular configuration, which implies 
 that some $x_i$ asymptotically reaches $l_j$ with $j \neq i,$ the limit of the respective
 correlation functions can be easily seen to coincide with the appropriate one-sided
 correlation function for the limit singular configuration, provided the convergence of $x_i$
 takes place on one side of $l_j$ only, otherwise the limit may fail to exist. The discontinuity
 may therefore arise because the one-sided correlation functions may differ on different
 sides.

 By Theorem 4 in \cite{SC08} it follows that for all 
 $(\bar l,\bar x) \in \Delta^g[k],\; k=1,2,$ we have 
 \begin{equation}\label{NISKIEK}
  \phi(\bar l,\bar x) = 1.
 \end{equation}
 The same paper gives general conditions for this relation
 to hold for $k > 2.$ Here we are interested in the general set-up for $k>2$ where
 it often happens that $\phi(\bar l,\bar x) \neq 1.$

% Dotad 18 IV 2009

 \paragraph{Representation theorem for edge correlation functions}
  The following theorem is the main result of this paper.
  \begin{theorem}\label{PWREPR}
   For each edge marker collection $(\bar l,\bar x)$ with $\bar x \subset D$ we
   have
   $$ \phi^{\cal M}(\bar l,\bar x) = {\Bbb E}\crop({\cal W}^{\cal M}_D[(\bar l,\bar x)]). $$
  \end{theorem}   
   It is useful at this point to compare Theorem \ref{PWREPR} with the simple observation
   that whenever all $x_i$'s in $(\bar l,\bar x)$ belong to a small ball $B_2(x_0,r)$ for
   some $x_0 \in {\Bbb R}^2$ then 
   \begin{equation}\label{PrzyblPhi}
    \phi(\bar l,\bar x) = N(\bar l,\bar x) (1+O(r)),
   \end{equation}
   with $N(\bar l,\bar x)$ defined as in Subsection \ref{ICF},
   which readily follows by (\ref{GIBBS1},\ref{GIBBS2}) and the definitions (\ref{MDEFI},\ref{MDEFI2})
   of edge correlation functions, taking additionally into account the local boundedness
   of the density $m(\cdot)$ of the activity measure ${\cal M}$ with respect to the
   Haar-Lebesgue measure $\mu$ which ensures that the Boltzmann factors 
   $\exp(-L^{\cal M}(\cdot))$ are $1+O(r)$ on $B_2(x_0,r),$ see also (12) in \cite{SC08}.
   It is interesting to note at this point
   that, by the consistency of ${\cal A}^{\cal M},$ the same would hold if we defined
   $N(\bar x,\bar l)$ to be the cardinality of $\Gamma(\bar x,\bar l) \cap D$ for any
   convex open $D \supset \bar x$ and thus $N(\bar x,\bar l)$ does not depend on the
   field domain $D$ as long as $D \supset \bar x$ !
   To proceed, observe that under the same conditions our Theorem \ref{PWREPR} yields
   \begin{equation}\label{PrzyblPhi2}
    \phi(\bar l,\bar x) = \crop({\cal W}^0_{D^r}[(\bar l,\bar x)]) (1+O(r))
   \end{equation}
   with ${\cal W}^0_{D^r}[(\bar l,\bar x)]$ defined as in Subsection \ref{ICF} and
   where $D^r$ is an r-dependent bounded convex domain of diameter $O(r)$ containing
   $B_2(x_0,r).$  
   This is because $O(r)$ is the probability that {\it at least one} turning/branching event occurs
   in the course of the polygonal web generating dynamics {\bf (W)} confined to $D^r,$
   again in view of the local boundedness of the density $m(\cdot).$ To proceed we note that
   both $N(\bar l,\bar x)$ and $\crop({\cal W}^0_{D^r}[(\bar l,\bar x)])$ are, by their
   definitions, invariant with respect to non-singular affine transforms of $(\bar l,\bar x).$
   Consequently, upon an appropriate re-scaling we can take $r$ in (\ref{PrzyblPhi}) and
   (\ref{PrzyblPhi2}) arbitrarily small. This way, upon comparing (\ref{PrzyblPhi}) and
   (\ref{PrzyblPhi2}) we have established
   \begin{lemma}\label{WN}
    For each edge marker collection $(\bar l,\bar x)$ we have 
    $$ \crop({\cal W}^0_D[(\bar l,\bar x)]) = N(\bar l,\bar x) $$
    for each $D \supset \bar x.$ In particular, $\crop({\cal W}^0_D[(\bar l,\bar x)])$
    does not depend on $D \supset \bar x.$
   \end{lemma}
   Note at this point that the relation (\ref{PrzyblPhi2}) can be further extended to
   produce a small $r$ expansion, with the coefficient at $r^k$ corresponding to instances 
   of polygonal webs with exactly $k$ turns/branchings. We do not pursue this topic here
   though because we are not aware of any natural geometric interpretations for the higher
   order terms of this expansion in the style of Lemma \ref{WN}. 

 \paragraph{Relaxing the {\it stop if separated} rule}
  As has already been remarked above, in contrast to the remaining rather strict dynamic
  rules, the {\bf (W:StopIfSeparated)} rule can be somewhat relaxed without
  affecting the validity of our Theorem \ref{PWREPR}. This is made more specific in the remark below.
  \begin{remark}\label{Sepa}
   A direct inspection of the proof of our representation Theorem \ref{PWREPR} below
   shows that the result stays valid if the {\bf (W:StopIfSeparated)} rule, requiring
   the web branch growth to cease {\it immediately} when its tip separates from the remaining
   ones, get replaced some other rule where the growth is stopped {\it only on separation}
   but not necessarily {\it immediately at separation.} The only natural constraints are that
   \begin{itemize}
    \item At each time moment in the course of the graphical construction the decision
          on whether to stop the growth of a separated branch either depends
          deterministically on the present configuration of branches or at least
          it is independent of the future evolutions of branches given the current
          branch configuration.
    \item Each branch eventually dies before or at the moment $s$ of becoming tangential to
          the current domain boundary $\partial D_s.$ 
   \end{itemize}   
   The first condition precludes unwanted dependencies whereas the second one is indispensible
   for the technical correctness of our constructions (the above tangency point is the point
   where the time flow direction changes along a branch, and the growth only occurs forward
   in time in our constructions).

   Note that a particular simple example of a stopping rule satisfying the above conditions is to kill
   each branch exactly at the time $s$ when it becomes tangential to the current boundary 
   $\partial D_s.$  
  \end{remark}
  In fact, an even deeper analysis of our argument below shows that further relaxation
  of the considered rule are admissible. We do not discuss these details in this paper 
  though as they are of no use for our present purposes.
  %in some cases
  %branch growth in {\bf (Crop)} dynamics can be stopped already when it separates from a suitably chosen
  %subcollection of other branches 

\section{Proof of the representation Theorem \ref{PWREPR}}\label{PMT}
 The purpose of this section is to prove our main Theorem \ref{PWREPR}. Our argument
 splits into several parts and requires some additional concepts.

 \paragraph{Edge marker process}
 To proceed towards establishing our representation theorem for edge correlation
 functions, we shall introduce a Markovian {\it edge marker} process whose construction
 can be to some extent regarded as a {\it backwards version} of the dynamic representation
 discussed in Section \ref{GDR}. Roughly speaking, the dynamic representation
 involved an {\it explosion} of the field from a single point up to the entire domain
 $D,$ whereas the edge marker process represents an evolution {\it back in time}
 and thus an {\it implosion} of the marker configuration, eventually to reach
 one of possible null states. In fact, the edge marker process will be seen
 to encode the construction of the polygonal web, see (\ref{WEM}) below,
 providing an interpolation
 between the original marker configuration $(\bar l,\bar x)$ and the full
 polygonal web ${\cal W}[(\bar l,\bar x)],$ whence the backwards time flow
 direction. To make all this specific, take the increasing
 family $(D_t)_{t \in [0,1]}$ of convex compacts satisfying
 {\bf (D1-4)} as chosen in Section \ref{PWEB}. Next, 
 consider a continuous time branching {\it edge marker process} 
 $\Psi_s := \Psi_{s;D},\; s \in [0,1],$ taking its values in finite families
 of {\it signed} and possibly empty edge marker configurations, with generic
 notation 
 $$ \Psi_s = 
     \left\{ \eta^{(p)} :(\bar{l}^{(p)}(s),\bar{x}^{(p)}(s)) = 
        (l^{(p)}_i(s),x^{(p)}_i(s))_{i=1}^{k_p} \right\}_{p=1}^m $$
 with $m$ and $k_p$ allowed to depend on the time $s$ and
 where $\eta^{(p)} \in \{ +1, -1 \}.$ As the notation suggests, the signs
 $\eta^{(p)}$ are attributed once and for all to their respective marker
 configurations and do not evolve in time. In addition, we always
 require that
 \begin{description}
  \item{\bf (EM:DomainShrink)} for each $s \in [0,1]$ the marker points
   $x^{(p)}_i(s)$ are all contained in the set $D_{1-s},$
 \end{description}
 that is to say the domain of the process $\Psi_s$ shrinks over time
 as informally discussed above.
 Given the initial state $\Psi_0$ with all marker points contained in
 $\bar D = D_1$ and, in addtion, assumed not to contain three different
 edge marker lines meeting at one point, the process $\Psi_s$ is gouverned
 by the following Markovian dynamics {\bf (EM)}, clearly preserving the
 latter property in view of the absolute continuity ${\cal M} \ll \mu.$
 \begin{description}
  \item {\bf (EM:DiscardIfSeparated)} If at some time $s$ an edge marker $(l^{(p)}_i(s),x^{(p)}_i(s))$
        has the property that $l^{(p)}_i(s)$ does not hit the convex hull 
        $$ C^{(p)}_i(s) := \conv(\{ x^{(q)}_j(s),\; q=1,\ldots,m; j=1,\ldots,k_q \} \setminus 
             \{ x^{(p)}_j(s) \}) $$
        generated by all the remaining $x^{(q)}_j(s)$'s in all marker configurations $(\bar{l}^{(q)}(s),
        \bar{x}^{(q)}(s)),\;$ $q=1,\ldots,m,$ in which case we say that the marker
        $(l^{(p)}_i(s),x^{(p)}_i(s))$ {\it separates} from $\Psi_s,$
        % its configuration $(\bar{l}^{(p)}(s),\bar{x}^{(p)}(s)),$
        then remove the marker $(l^{(p)}_i(s),x^{(p)}_i(s))$ from its configuration 
        $(\bar{l}^{(p)}(s),\bar{x}^{(p)}(s)).$
        
        If the removal of $(l^{(p)}_i(s),x^{(p)}_i(s))$ makes some other markers separate from $\Psi_s,$
        % their configuration,
        the {\bf (EM:DiscardIfSeparated)} rule applies for them as well and
        they are subsequently removed.  
  \item {\bf (EM:FoldInwards)} Between the critical moments listed below, at each time $s \geq 0$
        each edge marker point $x^{(p)}_i(s)$ lying at the boundary $\partial D_s$ \--- in
        which case we say the marker is in {\it boundary position} \--- moves along
        its corresponding marker line $l^{(p)}_i(s)$ so as to {\it always stay} at the boundary
        of the shrinking domain, that is to say $x^{(p)}_i(s+ds)$ arises as the intersection
        of $l^{(p)}_i(s)$ with $\partial D_{s+ds}.$ Note that marker points not lying at
        the boundary $\partial D_s$ {\it do not move} until they are eventually met by
        the boundary at some later time. 
  \item {\bf (EM:DiscardOnCollision)} If at some time $s$ two edge marker points $x^{(p)}_i(s)$
        and $x^{(p)}_j(s)$ within the same marker collection $(\bar{l}^{(p)}(s),\bar{x}^{(p)}(s))$
        collide (meet) along non-colinear directions $l^{(p)}_i(s)$ and $l^{(p)}_j(s),$
        the two markers are removed from the configuration.

        As a result of the above {\bf (EM:DiscardOnCollision)} rule some other edge markers
        may separate from $\Psi_s$ \--- in such a case the {\bf (EM:DiscardIfSeparated)} rule
        applies immediately and the markers get discarded.
  \item {\bf (EM:Kill)} In the course of the time interval $[s,s+ds],$
        an edge marker $(l^{(p)}_i(s),x^{(p)}_i(s))$ in boundary position, moving along   
        the segment $x^{(p)}_i[s,s+ds] :=  \overline{x^{(p)}_i(s) x^{(p)}_i(s+ds)},$
        gets removed from its configuration $(\bar{l}^{(p)}(s),\bar{x}^{(p)}(s))$
        with probability ${\cal M}([[x^{(p)}_i[s,s+ds]]]).$
        These updates are performed independently for all {\it different} boundary 
        edge markers throughout all marker configurations constituting the process
        $\Psi_s,$ yet they are performed {\it simultaneously} for all {\it equal} edge
        markers contained in different configurations, that is to say if
        $(x^{(p)}_i(s),l^{(p)}_i(s)) = (x^{(q)}_j(s),l^{(q)}_j(s))$ for some
        $p \neq q$ then the kill events during $[s,s+ds]$ coincide
        for both these markers. In other words, the killing mechanism
        is a.s. identical for all instances of an edge marker present in
        different configurations constituting $\Psi_s.$ 

        As a result of the above {\bf (EM:Kill)} rule some other edge markers may separate
        from $\Psi_s$ \--- in such a case the {\bf (EM:DiscardIfSeparated)} rule
        applies immediately and the markers get discarded.
  \item {\bf (EM:TurnAndBranch)} In the course of the time interval $[s,s+ds],$
        for each boundary edge marker $(l^{(p)}_i(s),x^{(p)}_i(s))$ moving along   
        the corresponding segment $x^{(p)}_i[s,s+ds] := \overline{x^{(p)}_i(s) x^{(p)}_i(s+ds)},$
        with probability ${\cal M}(dl)$ for $l \in [[x^{(p)}_i[s,s+ds]]]$ a 
        {\it turn-and-branch} update occurs in the direction of $l,$ which results
        in replacing the original marker configuration $(\bar{l}^{(p)},\bar{x}^{(p)})$
        by three {\it offspring} marker configurations in $\Psi_s,$ which are: 
        \begin{itemize}
         \item $\eta^{(p)}:(\bar{l}^{(p)}(s),\bar{x}^{(p)}(s))$ (unmodified offspring),
         \item $\eta^{(p)}:(\bar{l}^{(p)}(s),\bar{x}^{(p)}(s)) \setminus \{ (l^{(p)}_i(s),
                x^{(p)}_i(s)) \} \cup  \{ (l^{(p)}_i(s+ds) := l, x^{(p)}_i(s+ds)) \}$
                (directional update offspring \--- the original marker line $l^{(p)}_i(s)$
                {\it turns} in the direction of $l,$ that is to say $l^{(p)}_i(s+ds) := l$), 
         \item $-\eta^{(p)}:(\bar{l}^{(p)}(s),\bar{x}^{(p)}(s)) \cup \{ (l,
                x^{(p)}_i(s+ds)) \}$ (branched offspring \--- both the original marker and
                its directional update are present).
        \end{itemize}
        As in the case of {\bf (EM:Kill)} above, these
        updates are performed independently for all {\it different} boundary 
        edge markers throughout all marker configurations constituting the process
        $\Psi_s,$ yet they are performed {\it simultaneously} for all {\it equal} edge
        markers contained in different configurations, that is to say if
        $(x^{(p)}_i(s),l^{(p)}_i(s)) = (x^{(q)}_j(s),l^{(q)}_j(s))$ for some
        $p \neq q$ then the turning/branching updates during $[s,s+ds]$ coincide
        for both these markers. In other words, the turning/branching mechanism
        is a.s. identical for all instances of an edge marker present in
        different configurations constituting $\Psi_s.$ 
       
        Whenever the above {\bf (EM:TurnAndBranch)} update is performed with 
        $l \cap C^{(p)}_i(s+ds) = \emptyset,$ which results in the directionally
        updated edge marker separating from $\Psi_s,$ the rule
        {\bf (EM:DiscardIfSeparated)} applies immediately to the corresponding
        directional update and branched offsprings. 
  \item {\bf (EM:ForcedTurnAndBranch)} If a boundary edge marker point $x^{(p)}_i(s)$ crosses
        some edge marker line $l^{(p)}_j(s),\; j \neq i,\; l^{(p)}_j(s) \neq l^{(p)}_i(s),$
        during the time interval $[s,s+ds],$ a {\it forced turn-and-branch} update occurs in the
        direction of $l^{(p)}_j(s),$ which results
        in replacing the original marker configuration $(\bar{l}^{(p)},\bar{x}^{(p)})$
        by three {\it offspring} marker configurations in $\Psi_s,$ which are: 
        \begin{itemize}
         \item $\eta^{(p)}:(\bar{l}^{(p)}(s),\bar{x}^{(p)}(s))$ (unmodified offspring),
         \item $\eta^{(p)}:(\bar{l}^{(p)}(s),\bar{x}^{(p)}(s)) \setminus \{ (l^{(p)}_i(s),
                x^{(p)}_i(s)) \} \cup  \{ (l^{(p)}_i(s+ds) := l^{(p)}_j(s), x^{(p)}_i(s+ds)) \}$
                (directional update offspring \--- the original marker line $l^{(p)}_i(s)$
                {\it turns} in the direction of $l^{(p)}_j(s),$ that is to say $l^{(p)}_i(s+ds) 
                := l^{(p)}_j(s)$), 
         \item $-\eta^{(p)}:(\bar{l}^{(p)}(s),\bar{x}^{(p)}(s)) \cup \{ (l^{(p)}_j(s),
                x^{(p)}_i(s+ds)) \}$ (branched offspring \--- both the original marker and
                its directional update are present).
        \end{itemize}
        It should be noted at this point that,
        unlike in the usual {\bf (EM:TurnAndBranch)} discussed above,
        here {\bf (EM:DiscardIfSeparated)}
        has no chance of becoming applicable directly upon the update because
        along $l^{(p)}_j(s)$ there always exists a direction pointing at $C^{(p)}_i(s+ds),$
        namely that towards $x^{(p)}_j(s+ds).$ Moreover, observe also that, in effect of the
        so-defined {\bf (EM:ForcedTurnAndBranch)}
        update, $x^{(p)}_i(s+ds)$ and $x^{(p)}_j(s+ds)$ become {\it coupled} both in the
        directional update and branched offsprings. 
  \item {\bf (EM:UnbreakableCouplings)} If at some time $s$ in the course of their evolution
        two marker points
        $x^{(p)}_i(s)$ and $x^{(p)}_j(s)$ are coupled in their configuration then whenever they cease
        to be so in the original configuration or any of its offspring configurations, the 
        coupling-breaker configuration is instantly removed from $\Psi_s.$ 
        Note that this is equivalent to the rejection of configurations where
        \begin{itemize}
         \item a coupled edge marker gets killed in a collision
               {\bf (EM:DiscardOnCollision)} or in {\bf (EM:Kill)},
         \item a coupled edge marker modifies its direction in directional update offsprings
               arising in {\bf (EM:TurnAndBranch,ForcedTurnAndBranch)}. 
        \end{itemize}
        The unmodified and branched offsprings do not break couplings and
        neither can a coupling be broken in {\bf (EM:DiscardIfSeparated)} because
        coupled markers are never separated since they always point at their
        pair. 
        %as well
        %as configurations resulting from coupled edge markers upon their 
        %{\bf (EM:TurnAndBranch,ForcedTurnAndBranch)} directional updates
        %and their reduced offspring in {\bf (EM:TurnAndBranch)} \--- augmented
        %offsprings of such markers in %{\bf (EM:TurnAndBranch,ForcedTurnAndBranch)}
        %do not break couplings.  
        Note that in contrast to {\bf (EM:DiscardIfSeparated)} rule,
        where we discard individual edge markers, here we remove entire configurations.
 
        If at some time moment a marker point $x^{(p)}_i(s)$ reaches its coupled $x^{(p)}_j(s),$
        both markers coalesce and evolve henceforth as one, in particular all coupling restrictions
        cease to apply. Such meeting and coalescence may occur at the tangency point of the
        respective directional line to the current domain $\partial D_s$ in which case 
        the resulting single marker is instantly discarded from the system by application
        of {\bf (EM:DiscardIfSeparated)}.
\end{description}

  The above construction of the edge marker process may seem rather bizarre at the
  first look, but this is in fact a rather simple object. The moving boundary of
  the shrinking domain $D_s$ drives inwards polygonal branching random walks of
  constituent edge markers. The directional updating and branching mechanisms
  of these walks are determined by the activity measure ${\cal M}.$ The directions
  of the walks are always chosen to point at the convex hull generated by the remaining
  marker points in the process. If such a choice becomes impossible due to edge
  marker separation, the marker is discarded. Colliding edge markers are also
  discarded. One further rule is unbreakability of once established marker couplings,
  which is ensured by  rejecting coupling-breaker configurations. It is important to
  note that with probability $1,$ in the course of the dynamics all markers eventually
  \begin{itemize}
   \item either separate and are discarded in {\bf (EM:DiscardIfSeparated)},
   \item or disappear in collisions {\bf (EM:DiscardOnCollision)},
   \item or are killed in {\bf (EM:Kill)},
   \item or finally they have their configurations annihilated due
         to coupling breaks, as an application of {\bf (EM:UnbreakableCouplings)}.
  \end{itemize}
  Thus, at time $1$ the process
  $\Psi_1$ consists a.s. of signed empty marker configurations. 

  Another crucial observation, readily verified by comparing the {\bf (EM)} and
  {\bf (W)} dynamics, is that, on the event $\{ \Psi_0 = \{ + 1 : (\bar l,\bar x) \} \},$
  the union of trajectories traced by the constituent marker points $x^{(p)}_i(s),\;
  s \in [0,1],$ of $\Psi$ coincides with the web ${\cal W}[(\bar l,\bar x)],$ that is to say
  \begin{equation}\label{WEM}
   {\cal W}^{\cal M}_D[(\bar l,\bar x)] = \bigcup_{s \in [0,1]} \;
    \bigcup_{(\bar l^{(p)}(s),\bar x^{(p)}(s))
   \in \Psi_s} \; \bigcup_{(x^{(p)}_i(s),l^{(p)}_i(s)) \in (\bar l^{(p)}(s),\bar x^{(p)}(s))}
   \{ x^{(p)}_i(s) \}.
  \end{equation}
  Moreover, again by the construction, (the history of) each marker configuration present
  in $\Psi_1$ bijectively corresponds to a complete, minimal and normal collection of branches of
  ${\cal W}[(\bar l,\bar x)],$ whereas the complete abnormal branch collections correspond
  to marker configurations rejected in {\bf (EM:UnbreakableCouplings)} (this latter correspondence also
  becomes a bijection as soon as the minimality of branch collections is assumed).
   
  Recall now our assumption made above stating that the killing, directional updating and branching
  mechanisms, while independent for {\it different} markers, do coincide for equal markers. This
  assumption is clearly crucial for (\ref{WEM}) above to hold, but could be easily lifted
  without affecting the validity of a significant part of  the theory presented below.
  In fact, these mechanisms can also be coupled in any other way as soon as the Markovian
  property of the dynamics is preserved. We do not discuss this issue here though as the
  imposed coupling seems to be the most natural one and leading to simplest formulations.
  %, another
  %natural option being performing all turn/branch updates independent also among multiple
  %instances of the same edge marker. 
  
  Some concern may be raised by the branching nature of the {\bf (EM)} evolution \---
  a natural question is whether no cardinality explosions occur for $\Psi_s.$ 
  This possibility is easily excluded though, as stated below.
  \begin{lemma}\label{OGRWZROSTU}
   For each bounded open convex set $D$ and initial condition $\Psi_0$ there
   exists $c[D;\Psi_0] < + \infty$ such that
   $$ \forall_{s \in [0,1]} \;\; {\Bbb E} \card(\Psi_s) \leq c[D;\Psi_0]. $$
  \end{lemma}  
   To see it use first the relation (\ref{WEM}) and the discussion following it to
   conclude that, for all $s \in [0,1],$ the expectation ${\Bbb E}\card(\Psi_s)$
   bounded by $2^{\mbox{ number of branches of } {\cal W}[(\bar l,\bar x)]},$
   which is the maximum number of possible branch collections. Now, the
   expectation of this number is finite because ${\cal W}[(\bar l,\bar x)]$ arises
   from a (critical) binary branching process evolving during a {\it finite time interval}.
   Note that the forced branchings do not cause trouble here because they 
   only allow to extend {\it already existing} lines born at time $s$ to (at the furthest)
   the opposite side of $D_{1-s},$ whereas {\it new lines} are only born due to
   the usual critical branching.
    
 % Dotad 19 IV 2009   

 \paragraph{Correlation process}
  Having constructed the edge marker branching process $\Psi_s,$ we are now going
  to compose it with the correlation function to obtain the {\it edge correlation process}
  $\Phi_s,\; s \in [0,1].$ We put
  \begin{equation}\label{CORPROC}
   \Phi_s := \sum_{(\bar{l}^{(p)}(s),\bar{x}^{(p)}(s)) \in \Psi_s} \eta^{(p)}(s) 
   \phi(\bar{l}^{(p)}(s),\bar{x}^{(p)}(s)),\;\; s \in [0,1],
  \end{equation}
  which means defining the correlation process to be the sum of correlation functions
  for all marker configurations in $\Psi_s$ taken with their corresponding signs.
  For formal correctness it is convenient to adopt at this point the convention 
  that whenever at some time $s$ the marker point $x_i^{(p)}(s)$ is in a singular
  position lying on some $l_j^{(p)}(s),\; j \neq i,$ then the correlation function
  $\phi(\bar l^{(p)}(s),\bar x^{(p)}(s))$
  is interpreted as the {\it one-sided} correlation function in which the marker line
  $l^{(p)}_i(s)$ is replaced by half-line $[l^{\leftarrow}]^{(p)}_i(s)$ indicating
  the direction where $x^{(p)}_i(s)$ came from just before hitting $l^{(p)}_j(s).$
  %\--- note that with probability $1$ this is the opposite direction to that
  % indicated by the current directional
  % half-line $l^{\rightarrow}_i(s).$ 
  In view of (\ref{GIBBS1},\ref{GIBBS2}) and
  since ${\cal M} \ll \mu,$ with probability $1$ this is equivalent to putting
  in such case $\phi(\bar l^{(p)}(s),\bar x^{(p)}(s)) := \lim_{u \to s-} 
  \phi(\bar l^{(p)}(u),\bar x^{(p)}(u))$ for each $s$ where a singularity is
  reached. 
  
  Assume now that $\Psi_0 = \{ +1:(\bar l,\bar x) \}$ and recall from the discussion
  following the definition of the {\bf (EM)} dynamics that by the time $1$ the marker process
  $\Psi$ reaches a terminal state consisting entirely of signed {\it empty} marker
  configurations.
  In view of (\ref{WEM}) and the discussion following it, each such empty marker 
  configuration has its history encoded by some complete normal (minimal) branch collection
  in ${\cal W}[(\bar l,\bar x)],$ whereas complete abnormal (minimal) branch collections
  correspond to marker configurations rejected in {\bf (EM:UnbreakableCouplings)}.
  Moreover, the sign $\eta$ assigned to each empty marker configuration in $\Psi_1$
  can be readily verified, by induction in the number 
  of branchings in the course of crop graph creation, to be
  $(-1)^{\mbox{ number of branchings }}.$ Observing that for a complete branch
  collection, the number of branchings is simply the difference between the
  number of branches and number of roots (the latter coinciding with the cardinality of
  the initial marker collection $(\bar l,\bar x)$) and recalling that $\phi(\emptyset) = 1$
  we finally conclude from (\ref{CORPROC}) and (\ref{CROPFN}) that
  \begin{equation}\label{PHI1}
   \Phi_0 = \phi(\bar l,\bar x),\;\;\Phi_1 = \crop({\cal W}[(\bar l, \bar x)]).
  \end{equation}
 
 \paragraph{Martingale property of the correlation process} 
  With the notation introduced above, we claim that the edge correlation process
  is actually a martingale.
  \begin{lemma}\label{GLOWNE}
   The correlation process $(\Phi_s)_{s \in [0,1]}$ is a martingale with respect to the
   filtration ${\cal F}_s$ generated by the marker process $\Psi_s.$
  \end{lemma}   
   In view of the relation (\ref{PHI1}) Lemma \ref{GLOWNE} means we have just constructed
   a martingale interpolating between $(\bar l,\bar x)$ and $\crop({\cal W}[(\bar l, \bar x)]).$
   This immediately implies the assertion of Theorem \ref{PWREPR} upon putting
   $\Psi_0 := \{ +1:(\bar l,\bar x) \}.$ Thus, it remains to establish the crucial Lemma \ref{GLOWNE}.

 \paragraph{Proof of Lemma \ref{GLOWNE}}
 In view of the Markovian nature of the edge marker process
 $\Psi_s,$ to prove Lemma \ref{GLOWNE} it is enough to
 establish the desired martingale property
 at $s=0.$ Moreover, for simplicity we present our argument for the initial
 value $\Psi_0$ of the marker process consisting of a single marker configuration
 $(\bar{l}^{(1)}(0),\bar{x}^{(1)}(0)) := (\bar{l},\bar{x}) = (l_i,x_i)_{i=1}^k,$
 whence the general argument for higher cardinality of the initial condition is
 readily obtained by straightforward repetition for all configurations in $\Psi_0.$
 We assume without loss of generality that $x_1$ lies at the boundary $\partial D_1,$
 for only boundary markers undergo evolution under {\bf (EM)} dynamics. The remaining
 marker points $x_i,\; i=2,\ldots, k$ may lie both on $\partial D_1$ and in the
 interior of $D_1.$

 Using Lemma 3 in \cite{SC08}
 and the definition (\ref{MDEFI}) of correlation functions we see that if the
 line $l_1$ does not hit the convex hull $\conv(\{x_2,\ldots,x_k\})$ then 
 \begin{equation}\label{CIECIE}
  \phi(l_1,x_1;l_2,x_2;\ldots;l_k,x_k)
  = \phi(l_2,x_2;\ldots;l_k,x_k)
 \end{equation}
 which justifies the {\bf (EM:DiscardIfSeparated)} rule. 
 Thus, below with no loss of generality we constrain ourselves to the case where
 the {\bf (EM:DiscardIfSeparated)} rule does not apply during the period $[0,ds]$
 of the {\bf (EM)} evolution. 

 To proceed with our argument,
 we shall use the generalised dynamic representation described in Section \ref{GDR},
 with the same increasing family of convex compacts $(D_t)_{t \in [0,1]}$ as that used in
 the construction of the edge marker process. We also recall that $D_1 = \bar D$ where
 $D$ is the field domain. It should be recalled at this point that the generic time
 parameter $t$ of the dynamic representation is related by $t = 1-s$ to the usual
 time parameter $s$ of the edge marker and correlation processes. As already signalled
 above, the idea below is to relate the dynamics of the edge marker process with
 the graphical representation under inverted (backward) time flow. 

 If there is another marker point $x_i,\; i \neq 1,$ in boundary position with the
 property that $l_1$ and $l_i$ meet in $D_1 \setminus D_{1-ds}$ then, putting
 $i=2$ for notational clarity, by (\ref{MDEFI}) and the 
 {\bf (GE:VertexBirth)} dynamic rule we have,
 \begin{equation}\label{CIECIE2}
  \phi(l_1,x_1;l_2,x_2;l_3,x_3;\ldots;l_k,x_k)
  = (1+o(1)) \phi(l_3,x_3;\ldots;l_k,x_k).
 \end{equation}
 Clearly, the above event corresponds to marker point collision under {\bf (EM)}
 dynamics and thus (\ref{CIECIE2}) justifies the {\bf (EM:DiscardOnCollision)} rule.
 Keeping this in mind, below we only consider the case where $x_1$ does not
 collide with other marker points during $[0,ds].$
 
 To proceed, write the correlation-defining event
 $$ {\cal E}[dl_1,x_1;\ldots;dl_k,x_k] := \{ \forall_{i=1}^k 
     \exists_{e \in \Edges({\cal A}^{\cal M})} \; x_i \in e,\; l[e] \in dl_i \} $$
 as the intersection of ${\cal E}[dl_1,x_1]$ and ${\cal E}[dl_2,x_2;\ldots;dl_k,x_k]$
 where
 $$ {\cal E}[dl_1,x_1] := \{ 
     \exists_{e \in \Edges({\cal A}^{\cal M})} \; x_1 \in e,\; l[e] \in dl_1 \} $$
 and
 $$ {\cal E}[dl_2,x_2;\ldots;dl_k,x_k] := \{ \forall_{i=2}^k 
     \exists_{e \in \Edges({\cal A}^{\cal M})} \; x_i \in e,\; l[e] \in dl_i \}. $$
 We shall also denote by $x_0$ the intersection point of the marker line $l_1$ with
 $\partial D_{1-ds}.$ 
 With this notation and taking into account that ${\cal A}^{\cal M}_D = 
 {\cal A}^{\cal M} \cap D$ arises in the dynamic construction with $(D_t)_{t \in [0,1]}$
 as discussed above, we are now in a position to 
 represent ${\cal E}[dl_1,x_1;\ldots;dl_k,x_k]$ as the union of the following
 events, disjoint modulo a set of negligible probability, whose names were
 chosen to represent what happens if we
 {\it move along the field edge covering the marker} $(l_1,x_1)$ under the {\bf (EM)}
 dynamics inwards $D_1$ and towards $x_0.$
 \begin{description}
  \item{\bf (E:GoStraight)} 
        ${\cal E}[dl_2,x_2;\ldots;dl_k,x_k]$ occurs and an edge $e_1$ along $l_1$
        covering both $x_1$ and $x_0$ is present in the field ${\cal A}^{\cal M}_D.$
  \item{\bf (E:TurnOutwards)}
        ${\cal E}[dl_2,x_2;\ldots;dl_k,x_k]$ occurs, an edge $e_1$ along $l_1$ covering 
        $x_1$ is present but it does not reach $x_0,$ instead it turns at some point
        $x'$ of $\overline{x_1x_0}$ into another field edge $e'$ along a line $l'$ and in
        the direction consistent with that determined by the growth of $(D_t),$ that is to
        say $x'$ is the first point of $e'$ to be revealed by the growing window $(D_t).$ 
        Note that such an edge $e'$ can have only infinitesimally small length within
        $D_1 = \bar D$ since when moving from $x'$ in the direction indicated by the
        growth of $(D_t)$ we almost immediately encounter the boundary $\partial D.$
        Thus, we can say that $e'$ {\it points outside} the domain $\bar D$ and away from
        $D_{1-ds},$ whence the term {\it outward turn}. Often in such situations we shall
        also say that $e'$ extends {\it outwards} from $x'$ along $l'.$   

        In terms of the dynamic representation the occurrence of {\bf (E:TurnOutwards)}
        is equivalent to the occurrence of a {\bf (GE:VertexBirth)} vertex birth event
        at $x'$ between dynamic representation times $1-ds$ and $1,$ giving rise
        to the edges $e_1$ along $l_1$ and $e'$ along $l'.$  
  \item{\bf (E:TurnInwards)} 
        ${\cal E}[dl_2,x_2;\ldots;dl_k,x_k]$ occurs, an edge $e_1$ along $l_1$ covering 
        $x_1$ is present but it does not reach $x_0,$ instead it turns at some point
        $x'$ of $\overline{x_1x_0}$ into another field edge $e'$ along a line $l'$ and in the
        direction opposite to that determined by the growth of $(D_t),$ that is to say $x'$ is
        the last point of $e'$ to be revealed by the growing window $(D_t).$        
        In contrast to the above outward turn, here we turn in the opposite {\it inward}
        direction. In terms of the dynamic representation the occurrence of
        {\bf (E:TurnInwards)} is equivalent to the occurrence of a {\bf (GE:DirectionalUpdate)}
        at $x'$ where $e'$ extending along $l'$ turns into $e_1$ along $l_1.$

        For our
        considerations below it is convenient at this point to denote
        by $x''$ the intersection point of the inwards half-line ${e'}^{\rightarrow}$
        (starting at $x'$ and determined by $e'$)
        and $\partial D_{1-ds}.$ 
        %from $x'$ determined by $e',$ and such that $\dist(x_1,x')+\dist(x',x'')
        % = \delta$ (recall that $\delta = \dist(x_0,x_1)$).
        We say that $x''$ lies on $l'$ {\it inwards} from $x'$ in such cases. 
        Clearly, $x'' \in e'$ with probability $1-o(1).$  
 \end{description}
 %Dotad 30 IV 2009
 Observe now that, by the dynamic representation,
 \begin{itemize}
  \item  For {\bf (E:GoStraight)} we only consider the case that $x_1$ does not collide
         with any other marker point as the converse case has already been handled
         in (\ref{CIECIE2}). In this situation we have
         \begin{equation}\label{GStRnie}
            {\Bbb P}({\bf E:GoStraight}) =  (1+o(1))
            \left[ (1-{\cal M}([[\overline{x_1 x_0}]])) 
                   \sigma^{\cal M}[dl_1,x_0;dl_2,x_2;\ldots;dl_k,x_k]
                   \right. 
         \end{equation}
         $$ - \int_{l' \in [[\overline{x_1 x_0}]]}
              \sigma^{\cal M}[dl',x'';dl_1,x_0;dl_2,x_2;\ldots;dl_k,x_k] $$
         $$ - \sum_{l_j \in [[\overline{x_1 x_0}]],\; l_j \neq l_1}
              \left. \sigma^{\cal M}[dl_j,x''_j;dl_1,x_0;dl_2,x_2;\ldots;dl_k,x_k] \right]$$
         where $x''$ and $x''_j$ in the above integrals stand for the respective intersection
         points of $l'$ and $l_j$ with $\partial D_{1-ds},$ by definition lying inwards from
         $\{ x' \} := l' \cap \overline{x_0 x_1}$ and $\{ x'_j \} := l_j \cap \overline{x_0 x_1}.$ 
         To establish (\ref{GStRnie}) note first that on {\bf (E:GoStraight)} an edge $e_1$ 
         in the direction of $l_1$ is present at $x_0$ at the time $1-ds$ of the
         dynamic construction and, in the course of the dynamic construction,
         during the time period $[1-ds,1]$ the edge unfolds along $\overline{x_0 x_1}$
         eventually reaching $x_1.$ Consequently, there can be no {\bf (GE:DirectionalUpdate)}
         directional updates along $\overline{x_0 x_1}$ and neither can there be 
         {\bf (GE:Collision)} collisions with other already existing edges. However,
         according to the {\bf (GE)} dynamics such directional updates are possible
         along each line from $[[\overline{x_0 x_1}]],$ consequently the probability
         that none of these possible turns occurs yields the prefactor 
         $\exp(-{\cal M}([[\overline{x_1 x_0}]])) = (1+o(1))[1-{\cal M}([[\overline{x_1 x_0}]])]$
         whereas the remaining factor
         $\sigma^{\cal M}[dl_1,x_0;dl_2,x_2;\ldots;dl_k,x_k]$
         in the first term in the RHS of (\ref{GStRnie}) is the probability
         that the marker points $x_0,x_2,\ldots,x_k$ are covered by their corresponding
         edges as required. We claim that the resulting product
         \begin{equation}\label{PROOD}
          (1+o(1)) (1-{\cal M}([[\overline{x_1 x_0}]])) 
          \sigma^{\cal M}[dl_1,x_0;dl_2,x_2;\ldots;dl_k,x_k]
         \end{equation}
         represents the probability
         that the event ${\cal E}[dl_1,x_0;dl_2,x_2;\ldots;dl_k,x_k]$ holds and no 
         {\bf (GE:DirectionalUpdate)} turns occur along $\overline{x_0 x_1}.$ 
         To see it we observe that a possible directional update of this
         kind would yield, during the period $[1-ds,1]$ of the graphical
         construction, an outward edge $e'$ of infinitesimal length,
         almost immediately hitting the boundary $\partial D_1.$ Since
         for almost all time moments in the {\bf (EM)} dynamics the 
         distance between $x_1$ and other marker points $x_i,\; i > 1,$
         is strictly positive, during its short evolution under the 
         graphical construction dynamics the edge $e'$
         \begin{itemize}
          \item is overwhelmingly unlikely to have its birth event
                along the infinitesimal segment $\overline{x_0,x_1}$
                affected by the occurrence of 
                ${\cal E}[dl_1,x_0;dl_2,x_2;\ldots;dl_k,x_k],$
          \item has only a negligible chance of affecting the occurrence
                of the considered event ${\cal E}[dl_1,x_0;dl_2,x_2;\ldots;dl_k,x_k]$ because,
                in the course of the graphical construction, $e'$ is
                born {\it after} $x_0$ and other $x_i$'s contained in
                $D_{1-ds}$ get covered by the field and,
                in addition, $e'$ is, with overwhelming probability,
                too short to reach neighbourhood of any of the remaining
                $x_i$'s contained in $D_1 \setminus D_{1-ds}.$
         \end{itemize}       
         This nearly independence justifies taking the above product,
         as required. Next, we have to subtract the probability that 
         ${\cal E}[dl_1,x_0;dl_2,x_2;\ldots;dl_k,x_k]$ occurs and there
         are no directional updates {\bf (GE:DirectionalUpdate)} along $\overline{x_0 x_1}$
         but a {\bf (GE:Collision)} collision of the edge unfolding from $x_0$ along $l_1$
         with another already existing edge occurs on $\overline{x_0 x_1}.$
         There are two possible sources of such collisions
         \begin{itemize}
          \item $\overline{x_0 x_1}$ meets an inward edge $e'$ along some $l'$ non-colinear
                with any of the marker lines $l_j.$ Then the probability of the considered
                event is $$ (1+o(1)) \sigma^{\cal M}[dl',x'';dl_1,x_0;dl_2,x_2;\ldots;dl_k,x_k] $$
                with the prefactor $(1+o(1))$ due to the requirements that there be no 
                directional updates along $\overline{x_0 x_1}$ and that the inward
                edge $e'$ reaches $x'',$ which are negligibly unlikely to fail over the
                infinitesimal time interval $[0,ds].$ This expression
                corresponds to the second term in (\ref{GStRnie}) above.
          \item $\overline{x_0 x_1}$ meets an inward edge $e'$ along some marker line $l_j.$
                In analogy to the case above, here the probability of the considered
                event is 
                $$ (1+o(1)) \sigma^{\cal M}[dl_j,x''_j;dl_1,x_0;dl_2,x_2;\ldots;dl_k,x_k].$$
                This expression corresponds to the third term in (\ref{GStRnie}) above.
         \end{itemize} 
  \item If $x_1$ is not coupled with any other marker point $x_j$ then for {\bf (E:TurnOutwards)}
        we have   
        \begin{equation}\label{TOuRnie}
         {\Bbb P}({\bf E:TurnOutwards}) = (1+o(1)) {\cal M}(dl_1) {\cal M}([[\overline{x_1 x_0}]]) 
         \sigma^{\cal M}[dl_2,x_2;\ldots;dl_k,x_k].
        \end{equation}
         Indeed, the probability of the
         vertex birth {\bf (GE:VertexBirth)} at $x'$ as required for {\bf (E:TurnOutwards)}
         is ${\cal M}(dl_1) {\cal M}(dl'),$ with the notation as in the above
         definition of the event. Integrating over $l' \in [[\overline{x_1 x_0}]]$
         yields ${\cal M}(dl_1) {\cal M}([[\overline{x_1 x_0}]]).$ Moreover, 
         the occurrence of such a {\bf (GE:VertexBirth)} event at some dynamic
         construction time in $[1-ds,1],$ as yielding an infinitesimally short
         outward edge $e',$ is overwhelmingly unlikely to affect or be affected
         by the occurrence of ${\cal E}[dl_2,x_2;\ldots;dl_k,x_k]$  for precisely
         the same reasons as those justifying (\ref{PROOD}) above. 
         %as determined by the dynamic construction during its $[0,1-ds]$
         %time period.
         This nearly independence allows us to express 
         ${\Bbb P}({\bf E:TurnOutwards})$ as the product of 
         ${\cal M}(dl_1) {\cal M}(dl')$ and 
         $\sigma^{\cal M}[dl_2,x_2;\ldots;dl_k,x_k]$ with the extra
         $(1+o(1))$ coming also from the requirement that there be no {\bf (GE:DirectionalUpdate)}
         turns nor {\bf (GE:Collision)} collisions along $\overline{x' x_1}.$ 
        
         Note that the above conclusions are valid regardless of whether some marker line
         $l_j,\; j \neq 1,$ does cross $\overline{x_0 x_1}$ or not. Indeed, if
         such $l_j$ crosses $\overline{x_0 x_1}$ then the corresponding marker
         point $x_j$ lies in the inward direction from the intersection point
         for otherwise an {\bf (EM:DiscardOnCollision)} event would occur in
         the {\bf (EM)} dynamics, which we assumed not to be the case. Thus,
         a possible outward turn in the direction of $l_j,\; j \neq 1,$ occurring
         along $\overline{x_0 x_1}$ cannot yield an edge reaching and affecting the
         status of $x_j$ and the occurrence of ${\cal E}[dl_2,x_2;\ldots;dl_k,x_k].$
  \item  If $x_1$ is coupled with some other marker point
         $x_j,\; j \neq 1,$ along $l_1 = l_j$ then for {\bf (E:TurnOutwards)}
         we have ${\Bbb P}({\bf E:TurnOutwards}) = 0$ because the probability of obtaining
         in {\bf (GE:VertexBirth)} at $x'$ an edge $e_1$ {\it exactly colinear} with $l_j$ is
         zero since ${\cal M} \ll \mu.$           
  \item  If $x_1$ is not coupled with any other marker point $x_j$ then for {\bf (E:TurnInwards)}
         we have
         \begin{equation}\label{TInRnie}
          {\Bbb P}({\bf E:TurnInwards}) = (1+o(1)) {\cal M}(dl_1)
          \int_{l' \in [[\overline{x_1 x_0}]] } 
          \sigma^{\cal M}[dl',x'';dl_2,x_2;\ldots;dl_k,x_k].
         \end{equation}
         Note that if some $l_j,\; j \neq 1,$ crosses $\overline{x_0 x_1}$ then
         the integral in the RHS of (\ref{TInRnie}) above includes the singular 
         term ${\cal M}(dl_1) \sigma^{\cal M}[dl_j,x''_j;dl_2,x_2;\ldots;dl_j,x_j;\ldots;dl_k,x_k]$
         corresponding to the situation where $l' = l_j.$ To establish (\ref{TInRnie}) observe
         that the probability of our edge $e_1$ along $l_1$ arising due to a 
         {\bf (GE:DirectionalUpdate)} directional update at $x'$ on an inward edge $e'$ along $l'$ and of
         having ${\cal E}[dl_2,x_2;\ldots;dl_k,x_k]$ at the same time, is ${\cal M}(dl_1)$
         (directional update probability) times $\sigma^{\cal M}[dl',x'';dl_2,x_2;\ldots;dl_k,x_k]$
         (probability of ${\cal E}[dl_2,x_2;\ldots;dl_k,x_k]$ holding and of having an edge $e'$
          along $l'$ ending at $x'$ and thus passing through $x''$ modulo negligible measure set)
          times $(1+o(1))$ to take into account the requirement that there be no further
          {\bf (GE:DirectionalUpdate)} turns along $\overline{x' x_1}$ which is satisfied with
          overwhelming probability. In analogy to our previous considerations for (\ref{PROOD}),
          also here taking products of the above probabilities, modulo $(1+o(1)),$
          is well justified because
          the directional update at $x'$ has only a negligible chance of affecting the
          occurrence of  ${\cal E}[dl_2,x_2;\ldots;dl_k,x_k].$
          % which is determined in earlier stages
          %of the dynamic construction {\bf (GE)}.
          Integrating over $l'$ yields now
          (\ref{TInRnie}) as required.           
   \item  If $x_1$ is coupled with some other marker point $x_j,\; j \neq 1,$ along
          $l_1 = l_j$ then for {\bf (E:TurnInwards)} we have ${\Bbb P}({\bf E:TurnInwards}) = 0$
          because the probability of obtaining in {\bf (GE:DirectionalUpdate)} at $x'$ an
          edge $e_1$ {\it exactly colinear} with $l_j$ is zero since ${\cal M} \ll \mu.$
 \end{itemize}
 Putting now the above observations and formulae 
 (\ref{CIECIE},\ref{CIECIE2},\ref{GStRnie},\ref{TOuRnie},\ref{TInRnie})
 together, using the definition of the edge correlations (\ref{CORR})
 and recalling that $x_1$ is in boundary position as assumed, we see that,
 with the notation introduced in the above discussion
 \begin{itemize}
  \item If the marker $(l_1,x_1)$ separates from $(\bar l,\bar x)$ then 
        $$ \phi(l_1,x_1;l_2,x_2;\ldots;l_k,x_k) = \phi(l_2,x_2;\ldots;l_k,x_k). $$
  \item If there is another marker point $x_i,\; i \neq 1,$ say $i=2,$ in boundary position
        and with the property that $l_1$ and $l_2$ meet in $D_1 \setminus D_{1-ds}$ then
        $$ \phi(l_1,x_1;l_2,x_2;l_3,x_3;\ldots;l_k,x_k) = (1+o(1)) \phi(l_3,x_3;\ldots;l_k,x_k). $$
  \item Otherwise:
   \begin{itemize}
    \item If $x_1$ is not coupled with any other marker point $x_j,\; j \neq 1,$
  $$ 
  \sigma^{\cal M}(dl_1,x_1;\ldots;dl_k,x_k) = (1+o(1)) \left[ 
   (1 - {\cal M}([[\overline{x_1 x_0}]])) \sigma^{\cal M}[dl_1,x_0;\ldots;dl_k,x_k] \right. + $$ 
  $$ {\cal M}(dl_1) {\cal M}([[\overline{x_1 x_0}]]) \sigma^{\cal M}[dl_2,x_2;\ldots;dl_k,x_k] + $$
  $$ \int_{ l' \in [[\overline{x_1 x_0}]]} 
    [{\cal M}(dl_1) \sigma^{\cal M}[dl',x'';dl_2,x_2;\ldots;dl_k,x_k] - $$ $$ 
     \sigma^{\cal M}[dl',x'';dl_1,x_0;dl_2,x_2;\ldots;dl_k,x_k]] $$
  $$ - \sum_{j \neq 1,\;l_j \in [[\overline{x_0 x_1}]]} 
      \left. \sigma^{\cal M}[dl_j,x_j'';dl_1,x_0;dl_2,x_2;\ldots;dl_k,x_k] \right] $$
  whence, upon recalling the definition (\ref{MDEFI}) of the edge correlation function
  $\phi(\cdot),$
  $$
    \phi(l_1,x_1;\ldots;l_k,x_k) = (1+o(1)) \left[ 
     (1 - {\cal M}([[\overline{x_1 x_0}]])) \phi(l_1,x_0;\ldots;l_k,x_k) \right. + 
  $$
  $$ {\cal M}([[\overline{x_1 x_0}]]) \phi(l_2,x_2;\ldots;l_k,x_k) + $$
  $$ \int_{ l' \in [[\overline{x_1 x_0}]]} 
    [ \phi(l',x'';l_2,x_2;\ldots;l_k,x_k) - 
      \phi(l',x'';l_1,x_0;l_2,x_2;\ldots;l_k,x_k)] {\cal M}(dl') + $$
  $$ \sum_{j \neq 1,\;l_j \in [[\overline{x_0 x_1}]]}
    [\left. \phi(l_j,x_j'';l_2,x_2;\ldots;l_k,x_k) -
      \phi(l_j,x_j'';l_1,x_0;l_2,x_2;\ldots;l_k,x_k) ] \right] $$
  where the extra positive term in the last sum comes from separate treatment of 
  the case $l' = l_j$ in (\ref{TInRnie}), see the discussion directly following
  this display. This can be further rewritten as
  \begin{equation}\label{PHIzwykle}
    \phi(l_1,x_1;\ldots;l_k,x_k) = (1+o(1)) \left[ 
     (1 - 2 {\cal M}([[\overline{x_1 x_0}]])) \phi(l_1,x_0;\ldots;l_k,x_k) \right. + 
  \end{equation}
   $$ {\cal M}([[\overline{x_1 x_0}]]) \phi(l_2,x_2;\ldots;k_k,x_k) + $$
   $$ \int_{ l' \in [[\overline{x_1 x_0}]]} 
    [ \phi(l_1,x_0;\ldots) + \phi(l',x'';l_2,x_2;\ldots) - 
      \phi(l',x'';l_1,x_0;l_2,x_2;\ldots)] {\cal M}(dl') + $$
   $$ \sum_{j \neq 1,\;l_j \in [[\overline{x_0 x_1}]]}
    [\left. \phi(l_j,x_j'';l_2,x_2;\ldots;l_k,x_k) -
      \phi(l_j,x_j'';l_1,x_0;l_2,x_2;\ldots;l_k,x_k) ] \right]. $$
 \item Likewise, if $x_1$ is coupled with some other marker point $x_j,\; j \neq 1,$ then
  \begin{equation}\label{PHIsparowane}
   \phi(l_1,x_1;\ldots;l_k,x_k) = (1+o(1)) \left[ 
     (1 - {\cal M}([[\overline{x_1 x_0}]])) \phi(l_1,x_0;\ldots;l_k,x_k) \right. -
  \end{equation}
  $$  \int_{ l' \in [[\overline{x_1 x_0}]]}
      \phi(l',x'';l_1,x_0;l_2,x_2;\ldots;l_k,x_k) {\cal M}(dl') - $$
  $$  \sum_{l_i \in [[\overline{x_1 x_0}]],\; l_i \neq l_1} 
      \phi(l_i,x_i'';l_1,x_0;l_2,x_2;\ldots;l_k,x_k). $$
 \end{itemize}
 \end{itemize}
 To complete the above discussion we recall that in the complementary case
 $x_1 \not\in \partial D_1$ the point marker $x_1$ would not evolve under the
 {\bf (EM)} dynamics.

 To proceed we combine (\ref{CIECIE},\ref{CIECIE2},\ref{PHIzwykle},\ref{PHIsparowane})
 and extend these relations for all the other boundary marker points $x_i$ in $(\bar l,\bar x).$
 Recalling the evolution rules {\bf (EM)} for the edge marker process, taking into account
 that $\Psi_0 = \{ (\bar l,\bar x) \}$ and getting rid of
 the $(1+o(1))$ prefactors as $ds \to 0,$ we finally obtain the relation
 \begin{equation}\label{MARTdow}
  \Phi_0 = \phi(\bar l,\bar x) = {\Bbb E} \sum_{(\bar l^{(p)}(ds),\bar x^{(p)}(ds)) \in 
   \Psi_{ds}} \eta^{(p)}(ds) \phi(\bar l^{(p)}(ds),\bar x^{(p)}(ds)) = {\Bbb E}\Phi_{ds}.
 \end{equation}
 Observe in this context that 
 \begin{itemize}
  \item the second line in equation (\ref{PHIzwykle}) corresponds
        to {\bf (EM:Kill)} rule,
  \item the third line there to {\bf (EM:TurnAndBranch)} rule,
  \item and the fourth line to {\bf (EM:ForcedTurnAndBranch)}.
 \end{itemize}
 Likewise, the absence of 
 certain terms in the coupled version (\ref{PHIsparowane}) of eq. (\ref{PHIzwykle})
 corresponds to annihilation of coupling-breaker configurations in 
 {\bf (EM:UnbreakableCouplings)}.
 Clear\-ly, the crucial relation (\ref{MARTdow})
 admits straightforward extensions for more general initial
 conditions and all time moments between $0$ and $1,$ as discussed at the beginning
 of our proof. Thus, (\ref{MARTdow}) implies in particular that 
 $$ \tilde{\Phi}_s := \sum_{(\bar{l}^{(p)}(s),\bar{x}^{(p)}(s)) \in \Psi_s} 
    \phi(\bar{l}^{(p)}(s),\bar{x}^{(p)}(s)),\;\; s \in [0,1],
 $$  
 where all $\eta^{(p)}$ signs are converted into pluses, is a positive submartingale.
 Using Lemma \ref{OGRWZROSTU} we see it is a uniformly integrable submartingale. 
 Therefore, noting that $|\Phi_s| \leq \tilde{\Phi}_s,$ we conclude from (\ref{MARTdow}) 
 that $(\Phi_s)_{s \in [0,1]}$ is a martingale, which completes the proof of Lemma \ref{GLOWNE}. $\Box$

 \paragraph{Completing the proof of Theorem \ref{PWREPR}}
  With Lemma \ref{GLOWNE} established, we are now in a position to use the martingale
  representation combined with the relation (\ref{PHI1}) to complete the proof of Theorem
  \ref{PWREPR} as discussed next to the statement of Lemma \ref{GLOWNE} above. $\Box$

\section{Proof of Lemma \ref{NCOR}}\label{PNCOR}
 Our proof is based on the so-called defective disagreement loop dynamics developed
 in Subsection 6.4 of \cite{SC08}. Since Lemma \ref{NCOR} is of a purely 
 technical rather than conceptual nature and the quite complicated defective
 diagreement loop dynamics finds no further applications in this paper,
 we decided not to present its several pages long details here, referring
 the reader to results of Subsection 6.4 and Section 7 in \cite{SC08}
 instead. To proceed, with $(\bar l, \bar x) = (l_1,x_1,\ldots,l_k,x_k)$
 write
 $$ {\cal E}[dl_1,x_1;\ldots;dl_k,x_k] := \{ \forall_{i=1}^k 
     \exists_{e \in \Edges({\cal A}^{\cal M})} \; x_i \in e,\; l[e] \in dl_i \} $$
 for the correlation-defining events. We claim that for
 $(l_i,x_i)_{i=1}^k \in \Delta^g[k]$ in general position, upon fixing
 $l_1,x_1;\ldots ;l_{k-1},x_{k-1},$ the function 
 \begin{equation}\label{FKind}
   (l_k,x_k) \mapsto {\Bbb P}({\cal E}[dl_k,x_k]|{\cal E}[dl_1,x_1;\ldots;dl_{k-1},x_{k-1}])
    / {\cal M}(dl_k)
 \end{equation}
 is well defined, continuous and consequently locally bounded. Clearly, this will
 imply the statement of Lemma \ref{NCOR} upon inductive application 
 for $(\bar l, \bar x) \in \Delta^g[k].$
 The existence and local boundedness for $(\bar l, \bar x) \in \Delta^s[k]$
 will then follow as well by noting that the edge correlation for a singular
 configuration is bounded above by the sum of all corresponding one-sided
 correlations and by repeating the existence and local boundedness argument
 given below for the case of one-sided correlations, which is a straightforward
 repetition omitted here to avoid unnecessary technicalities. 
 
 To establish our claim for (\ref{FKind}) we use the defective disagreement 
 loop dynamics of Subsection 6.4 in \cite{SC08} with directional updating
 principle induced by an arbitrary growing family $(D_t)_{t \in [0,1]}$ as
 in {\bf (D1-4)} satisfying in addition $D_0 = \{ x_1 \}$ so that the anchor
 point ${\Bbb A}(l_1)$ coincides with $x_1.$ In analogy to the proof of 
 Theorem 4 in \cite{SC08}, see also Theorem 10 there, the conditional
 law of the polygonal field ${\cal A}^{\cal M}$ on the event
 ${\cal E}[dl_1,x_1;\ldots;dl_{k-1},x_{k-1}]$
 is invariant with respect to the following reversible conditional version of
 the defective disagreement loop dynamics, with $s$ standing for the
 corresponding time parameter in which the dynamics unfolds: 
 \begin{description}
  \item{\bf (Create)} With intensity ${\cal M}(dl_k) ds,$ on the event
       $\neg {\cal E}[dl_k,x_k]$ (that is to say if there is no field 
       edge along $dl_k$ containing $x_k$) attempt to emit from $x_k$
       a disagreement path with initial creation phase directed along $l_k.$
       Should the so generated path result in a configuration violating
       ${\cal E}[dl_1,x_1;\ldots;dl_{k-1},x_{k-1}],$ discard the update,
       otherwise accept it.
  \item{\bf (Annihilate)} With intensity $ds,$ on the event ${\cal E}[dl_k,x_k]$
       (that is to say if there is a field edge along $dl_k$ containing $x_k$) 
        attempt to emit from $x_k$ a disagreement path with initial annihilation
        phase directed along $l_k.$ Should the so
        generated path result in a configuration violating
        ${\cal E}[dl_1,x_1;\ldots;dl_{k-1},x_{k-1}],$ discard
        the update, otherwise accept it.   
 \end{description}
 Note that there are no update failures arising due to cycle formation along
 disagreement paths in this dynamics, because the chosen directional updating
 rule comes from a generalised dynamic construction, see Section 6.4 in 
 \cite{SC08} for details. Denote now by $\pi_{\rm create}$ the conditional probability of
 a succesful {\bf Create} update attempt during time interval $(s,s+ds)$ on
 the event $\neg {\cal E}[dl_k,x_k].$ Likewise, write $\pi_{\rm annihilate}$
 for the respective conditional probability for {\bf Annihilate} update on
 ${\cal E}[dl_k,x_k].$ 
 Clearly, by detailed balance for ${\cal E}[dl_k,x_k],$ we have
 \begin{equation}\label{DZIELENIE}
   {\Bbb P}({\cal E}[dl_k,x_k]|{\cal E}[dl_1,x_1;\ldots;dl_{k-1},x_{k-1}])
   = (1+o(1)) \pi_{\rm create} / \pi_{\rm annihilate}.
 \end{equation}
 Note that $\pi_{\rm create} \leq {\cal M}(dl_k) ds,$ and that
 $\pi_{\rm annihilate} \geq c ds$ for some $c$ uniformly positive
 with respect to small local displacements of $x_k$ and
 $l_k$ because there is some positive probability that the
 disagreement path initiated by annihilating the edge at $x_k$
 does not hit any $x_i,\; i < k,$ and thus does not lead to
 the violation of ${\cal E}[dl_1,x_1;\ldots;dl_{k-1},x_{k-1}].$  
 Consequently, we conclude from (\ref{DZIELENIE}) that the function in
 (\ref{FKind}) is well defined and locally bounded on $\Delta^g[k].$
 Its required continuity follows also by (\ref{DZIELENIE}) in view
 of the assumed continuity of the activity measure density
 $m = d{\cal M}/d\mu.$ We have thus established the desired
 properties of the conditional correlation in (\ref{FKind})
 which completes the proof of the lemma. $\Box$     

%     \leq k {\cal M}(dl_k). $$
% \end{lemma}
% Clearly, Lemma \ref{NCOR} follows from Lemma \ref{BDCOND} by straightforward induction.

\paragraph{Acknowledgements}
 The author acknowledges the support from the Polish Minister of Science and Higher 
 Education grant N N201 385234 (2008-2010).


\begin{thebibliography}{999}
  \bibitem{A1}
       {\sc Arak, T.} (1982) On Markovian random fields with finite number of values,
       {\it 4th USSR-Japan symposium on probability theory and mathematical statistics,
           Abstracts of Communications}, Tbilisi.
  \bibitem{AS1}
       {\sc Arak, T., Surgailis, D.} (1989) Markov Fields with Polygonal Realizations,
       {\it Probab. Th. Rel. Fields} {\bf 80}, 543-579.
  \bibitem{AS2}
       {\sc Arak, T., Surgailis, D.} (1991) Consistent polygonal fields,
       {\it Probab. Th. Rel. Fields} {\bf 89}, 319-346.
  \bibitem{ACS}
        {\sc Arak, T., Clifford, P., Surgailis, D.} (1993) Point-based polygonal models
        for random graphs, {\it Adv. Appl. Probab.} {\bf 25}, 348-372.
  \bibitem{AC07} {\sc Arias-Castro, E.} (2007) Interpolation of random hyperplanes,
        {\it Electronic Journal of Probability}, {\bf 12}, 1052-1071.
  \bibitem{CN94}
       {\sc Clifford, P., Nicholls, G.} (1994).
        A Metropolis sampler for polygonal image reconstruction.
        {\it available at}:\\
        {\tt http://www.stats.ox.ac.uk/~clifford/papers/met\_poly.html}
  \bibitem{FFG1}
      {\sc Fern\'andez, R., Ferrari, P., Garcia, N.} (1998)
        Measures on contour, polymer or animal models. A probabilistic
        approach. {\it Markov Processes and Related Fields} {\bf 4}, 479-497.
  \bibitem{FFG2}
      {\sc Fern\'andez, R., Ferrari, P., Garcia, N.} (2001)
       Loss network representation of Ising contours
       {\it Ann. Probab.} {\bf 29}, 902-937. 
  \bibitem{FFG3}
      {\sc Fern\'andez, R., Ferrari, P., Garcia, N.} (2002) Perfect simulation for
        interacting point processes, loss networks and Ising
        models. {\it Stoch. Proc. Appl.} {\bf 102}, 63-88.
  \bibitem{KLS05} {\sc Kluszczy\'nski, R., Lieshout, M.N.M.~van, Schreiber, T.}
       (2005) An algorithm for binary image segmentation using polygonal Markov fields.
       In: F.\ Roli and S.\ Vitulano (Eds.), Image Analysis and Processing, Proceedings
       of the 13th International Conference on Image Analysis and Processing.
       {\it Lecture Notes in Comput. Sci.} {\bf 3615}, 383-390. 
  \bibitem{KLS07} {\sc Kluszczy\'nski, R., Lieshout, M.N.M.~van, Schreiber, T.} 
       (2007) Image segmentation by polygonal Markov fields. {\it Ann. Inst. Statist. Math.},
       {\bf 59}, 465-486.
  \bibitem{LS07} {\sc Lieshout, M.N.M.~van, Schreiber, T.} (2007) 
       Perfect simulation for length-interacting polygonal Markov fields in
       the plane, {\it Scand. Journal of Statistics}, {\bf 34}, 615-625.
  \bibitem{LG} {\sc Liggett, T.} (1985).
       {\it Interacting particle systems}.
       Springer-Verlag, New York.
  \bibitem{MM} {\sc Malyshev, V.A., Minlos, R.A.} (1991) {\it Gibbs random fields: cluster expansions},
       Mathematics and Its Applications (Soviet Series), Kluwer.
  \bibitem{N1}
       {\sc Nicholls, G.K.} (2001) Spontaneous magnetization in the plane,
        {\it Journal of Statistical Physics}, {\bf 102}, 1229-1251.
  \bibitem{SC05}
       {\sc Schreiber, T.} (2005) Random dynamics and thermodynamic limits for 
        polygonal Markov fields in the plane, {\it Advances in Applied
        Probability} {\bf 37}, 884-907.
  \bibitem{SC06}
       {\sc Schreiber, T.} (2006) Dobrushin-Koteck\'y-Shlosman theorem for polygonal
        Markov fields in the plane, {\it Journal of Statistical Physics}, {\bf 123},
        631-684.
 \bibitem{SC08} {\sc Schreiber, T.} (2008) Non-homogeneous polygonal fields in the plane:
                graphical constructions and geometry of higher order correlations,
                {\it Journal of Statistical Physics}, {\bf 132}, 669-705. 
 \bibitem{SL08} {\sc Schreiber, T., Lieshout, M.N.M.~van}
       (2008) Disagreement loop and path creation/annihilation algorithms for
       binary planar Markov fields with applications to image segmentation,
       {\it submitted}
  \bibitem{SU1}
       {\sc Surgailis, D.} (1991) Thermodynamic limit of polygonal models, {\it Acta applicandae
        mathematicae}, {\bf 22}, 77-102.
 \end{thebibliography}
\end{document}